\documentclass[12pt]{amsart}
\usepackage{graphicx}
\usepackage{todonotes}

\usepackage{indentfirst, csquotes}

\usepackage[top=.8in, left=0.8in, bottom=.8in, right=0.8in]{geometry}

\newtheorem{theorem}{Theorem}[]

\newtheorem{lemma}[theorem]{Lemma}

\newtheorem{corollary}[theorem]{Corollary}

\newtheorem{remark}[theorem]{Remark}

\usepackage{amssymb,amsthm,amsmath}
\usepackage{xcolor,paralist,titlesec,fancyhdr,etoolbox}
\usepackage[colorlinks=true, citecolor=black, urlcolor=blue]{hyperref}
\usepackage{pdflscape} 

\usepackage{easybmat}

\titleformat{\section}[display]{\normalfont\huge\bfseries\centering}{}{10pt}{\Large}
\titlespacing*{\section}{0pt}{0ex}{0ex}

\hypersetup{ colorlinks=true, linkcolor=black, filecolor=black, urlcolor=black }

\newcommand{\R}{\mathbb{R}}
\newcommand{\C}{\mathbb{C}}

\newcommand{\la}{\lambda}
\newcommand{\g}{\mathfrak{g}}
\newcommand{\h}{\mathfrak{h}}

\newcommand{\p}{\mathfrak{p}}
\newcommand{\kk}{\mathfrak{k}}
\newcommand{\so}{\mathfrak{so}}

\newcommand{\q}{\mathfrak{q}}
\newcommand{\su}{\mathfrak{su}}
\newcommand{\un}{\mathfrak{u}}

\renewcommand{\sp}{\mathfrak{sp}}

\DeclareMathOperator*{\proj}{proj}

\begin{document}
	\title{ON THE PRESCRIBED RICCI CURVATURE OF NONCOMPACT HOMOGENEOUS SPACES WITH TWO ISOTROPY SUMMANDS} 
	\author[]{Dustin Gaskins}
	\date{\today}
	\address{Address}
	\email{example@mail.com}
	\maketitle
	
	\let\thefootnote\relax
	
	\begin{abstract}
		This work studies simply connected, noncompact $G/H$ in which $G$ is semi-simple, $H$ is connected, and $G/H$ has two irreducible summands. Here, we classify all such spaces and we provide solutions to the so-called Prescribed Ricci Curvature problem for all such spaces.
	\end{abstract} 
	
	\tableofcontents
	
	\bigskip

	\section{1. Introduction}\label{introduction}
	
	\vspace{.2in}
	
	\noindent Let $G/H$ being a simply connected homogeneous manifold with $H$ connected, a $G$-invariant metric $g$, and Ricci curvature $ric_g(.,.)$. The past century has contained a considerable amount of effort to understand homogeneous spaces and the geometry of said spaces. Much of this interest can be summed up by the question, \textit{what is the relationship between the geometry and the topology of homogeneous spaces?} One such relationship that has received a great deal of attention is the relationship between compactness and curvature. Some examples include Myer's theorem in the more general Riemannian geometry setting, stating that if a complete Riemannian manifold has Ricci curvature bounded below by a positive constant, then the space is compact (see \cite{doCarmo}). Years later in the 1970s, another result of this nature was provided by  B. Bergery in Theorem 2 of \cite{BB}, stating that a homogeneous space having universal cover diffeomorphic to Euclidean space is equivalent to being flat or having strictly negative scalar curvature. 
	\\
	\\Over the past several decades, many have been investigating the effects of having an Einstein metric (i.e $ric_g = \la g$ for some $g$) on a homogeneous space. Regarding the existence of Einstein metrics in the compact setting, one can refer to the survey \cite{CompEinSurv}. In the non-compact setting, though, a recent and exciting development of considerable magnitude has been the solution to the Alekseevskii conjecture (see \cite{AlekSol}), proving that a connected homogeneous Einstein space with negative scalar curvature is diffeomorphic to Euclidean space.
	\\
	\\Other recent developments include investigations (see \cite{PRPSurvey, Pul1, Pul5, ArrPRPC, ArrPRPN, LauPRP, Butt}) into finding solutions to the so-called \textit{Prescribed Ricci Curvature Problem} (PRP) which seeks to find the solutions  to Eqns. (\ref{PRP1}) and (\ref{PRP2}), with high emphasis on the latter.
	\\
	\begin{align}\label{PRP1}
		ric_g(.,.) &= T(.,.)
	\end{align}
	\begin{align}\label{PRP2}
		ric_g(.,.) &= cT(.,.)
		\\\notag
	\end{align}
	In Eqn. (\ref{PRP1}), the goal is to find all $T$ for which there is a metric $g$ such that $ric_g = T$, and in Eqn. (\ref{PRP2}) the goal is to find all $T$ for which there is a $c\in\R$ and a metric $g$ with $ric_g = cT$. Geometrically, Eqn. (\ref{PRP1}) can be thought of as finding the image of $ric$  and Eqn. (\ref{PRP2}) the image of $ric$ up to scaling. In the latter equation, when possible, it is also desirable to provide the $c$ value for a given solution $T$.
	\\
	\\While much work investigating the PRP has been done in the compact setting, little has been done in the noncompact setting (as the survey of developments in \cite{PRPSurvey} notes). 
	\\
	\\Similarly, but shifting perspectives slightly, in the past couple of decades there has been a considerable amount of effort to understand compact spaces in which $G/H$ has two irreducible isotropy summands. Indeed, with the isotropy irreducible spaces being settled in both the compact and noncompact cases (see \cite{JWolf}) and the geometry well understood, it is reasonable to consider the next easiest set of spaces to investigate: those with two irreducible summands. 
	\\
	\\An expansion in the direction of two irreducible summands was notably done in \cite{Kerr} in which a classification of such $G/H$ with $G$ simple is obtained (with corrections in Theorem A.1 of \cite{IsoCorect} and Remark 6.1 of \cite{LauretTwoCorrection}), as well as a list of which spaces contain Einstein metrics (i.e have metrics in which $ric = \la g$). In addition to this work of Dickinson's and Kerr's, there has also been a thorough investigation into the Ricci flow of compact $G/H$ with two inequivalent irreducible isotropy summands (\cite{BuzIso}), and more recently, a solution to Eqn. (\ref{PRP2}) for compact $G/H$ ($G$ not necessarily simple) with two inequivalent isotropy irreducible summands in \cite{Pul4}.
	\\
	\\With such rich, recent developments in the compact setting for both the PRP and the case of two irreducible summands, the present work seeks to investigate the corresponding noncompact setting. Working specifically in the setting in which $G/H$ is simply connected,  $G$ is semi-simple, and $H$ is connected, we first prove the following theorem:
	\begin{theorem}
		Let $G/H$ be simply connected with $G$ a connected semi-simple Lie group with no compact factors and $H \subset G$, a compact, connected subgroup. If $G/H$ has exactly two irreducible representations then $G/H$ is described by one of the following:
		\begin{enumerate}
			\item $G$ and $H$ have Lie algebras $\g = \g_1 \oplus \g_2$ and $\h = \kk_1 \oplus \kk_2$ where $\g_i$ is noncompact simple and $\kk_i$ is the maximal compact in $\g_i$. In this case, $G/H$ is a symmetric space. 
			\item $G$ has a simple Lie algebra $\g$, and $H$ has Lie algebra $\h \subsetneq \kk$ where $\kk$ is the maximal compact in $\g$. A classification of such $G/H$ is determined by the $(\g, \kk, \h)$ triple belonging to Tables \ref{TwoSummandsTable1} through \ref{TwoSummandsTable4}. In this case, $G/H$ is not a symmetric space.
			\\
			\\
		\end{enumerate}
	\end{theorem}
	
	\noindent After we prove this theorem, we then provide solutions to Eqns. (\ref{PRP1}) and (\ref{PRP2}). The $G$ semi-simple and not simple setting turns out to be quite simple, but the $G$ simple setting is not so simple. A brief and non-comprehensive comparison of solutions to Eqn. (\ref{PRP2}) for these two contexts makes this clear:
	\\
	\\
	\begin{theorem}
		If $\g$ is not simple then $ric = cT$ for $c > 0$ if and only if $T$ is a scalar multiple of $ric_{\langle . , .\rangle}$ where $\langle . , . \rangle$ is our fixed inner product, and $c = \frac{T(x,x)}{ric_{\langle . , . \rangle}(x,x)}$ for some $x \in \p = \p_1 \oplus \p_2$. Moreover, our $T < 0$ in this case.
	\end{theorem}
	
	\begin{theorem}
		Let $\g = \h \oplus \p_1 \oplus \p_2$ be noncompact simple with $[\p_1, \p_1] \not\subset \h$ (except $\g = \so(1,7)$). For $G/H$ in this case, the equation $ric = cT$ for $c> 0$ has a solution if and only if $(t_1, t_2)$ is a pair satisfying $$\frac{t_1}{t_2} \leq \frac{d_2^2}{d_1p_2} - \frac{d_2}{d_1}\sqrt{\frac{d_2^2}{p_2^2} + \dfrac{2d_1 - p_1 - 2p_2}{p_2}} \leq 0$$
		where $t_1 < 0$ and $t_2 > 0$. (Equality on the right only occurs when $\p_1$ is trivial. See Corollary \ref{ricequalscTtrivialp2} for that setting).
		\\
		\\When the above inequality is satisfied, there is always one solution which can be obtained by $c_+$ (where $c_+$ takes the $+$ in (\ref{cvaltop}) below) and $(x_1, x_2)$, the pair unique up to scaling given by $\frac{x_2}{x_1} = \frac{-2d_2}{p_2}(c_+ t_2 + \frac{1}{2}).$ 
		\\
		\\In addition to this one solution, there is a second solution if and only if our pair $(t_1, t_2)$ with $t_1 < 0$ and $t_2 > 0$ satisfies
		$$-\dfrac{2d_1 - p_1 - 2p_2}{2d_1} < \frac{t_1}{t_2} < \frac{d_2^2}{d_1p_2} - \frac{d_2}{d_1}\sqrt{\frac{d_2^2}{p_2^2} + \dfrac{2d_1 - p_1 - 2p_2}{p_2}}.$$
		The second solution is obtained by $c_-$ (where $c_-$ takes the $-$ in \ref{cvaltop} below) and $(x_1, x_2)$, the pair unique up to scaling given by $\frac{x_1}{x_2} = \frac{-2d_2}{p_2}(c_- t_2 + \frac{1}{2}).$
		\begin{align}\label{cvaltop}
			c = \dfrac{-\left(\frac{d_2^2}{d_1p_2}t_2 - t_1\right) \pm \sqrt{\left(\frac{d_2^2}{d_1p_2}t_2 - t_1\right)^2 - 4\left(\frac{d_2^2}{d_1p_2}t_2^2\right)\left(\frac{d_2^2}{4p_2d_1} + \dfrac{2d_1-p_1-2p_2}{4d_1}\right)}}{2\frac{d_2^2}{d_1p_2}t_2^2} 
		\end{align}
		Moreover, $c_+ = c_-$ when the discriminant in $c$ is zero, and this happens precisely when $$\frac{t_1}{t_2} = \frac{d_2^2}{d_1p_2} - \frac{d_2}{d_1}\sqrt{\frac{d_2^2}{p_2^2} + \dfrac{2d_1 - p_1 - 2p_2}{p_2}}.$$
	\end{theorem}
	
	\noindent In addition to working in the noncompact setting, there are some other distinctions between our work and prior works in the compact setting. One notable distinction is that we do not restrict ourselves to the assumption that our two irreducible summands are inequivalent. This turns out to be an exceptional case working with $SO_0(1,7)/G_2$, and the PRP results can be found in Theorems \ref{RiceqT} and \ref{RiceqcT}.
	\\
	\\Another notable distinction regards our methodology. In the aforementioned works, the dominant method for finding solutions to Eqn. (\ref{PRP2}) is done by utilizing a restriction of metrics and applying variational methods on those metrics. To be precise, in works such as \cite{Pul1, Pul4, ArrPRPC, ArrPRPN}, they work with those $G$-invariant metrics belonging to:
	\begin{align*}
		M_T = \{g: tr_g T = 1\}.
	\end{align*}
	In restricting themselves to such metrics, they utilize a variational approach provided in the following lemma:
	\begin{lemma}
		Given $g \in M_T$, there is a solution to Eqn. (\ref{PRP2}) for some $c\in \R$ if and only if $g$ is a critical point of the restriction of the scalar curvature function to $M_T$.
	\end{lemma}
	\noindent While there is considerable strength to this lemma, it also has some limitations. For a given solution $T_0$, there is a metric $g$ and $c = \pm S(g)$, where $S(g)$ is the scalar curvature. Unless, however, one's method for finding the critical points also reveals the metric in terms of $T_0$, this will leaves one wondering what $g$ to choose or what $ric_g$ is for the solution $T_0$. Thus, it would be more ideal to be able to say what $c$ is in terms of $T_0$ or to know what $g$ is such that $ric_g = cT_0$.
	\\
	\\In the present work, we do not use variational methods, and we are able to determine what $c$ is in terms of a given solution $T_0$. While our approach certainly has its limitations as well, it makes up for said limitations by providing an exhaustive solution to both Eqns. (\ref{PRP1}) and (\ref{PRP2}) for $G/H$ in which $G$ is semi-simple, $H$ is connected, and $G/H$ has two irreducible summands, with the $c$ values determined as well. Our methods also describe the metric in terms of $c$ and $T$, except in the special case of $SO_0(1,7)/G_2$ in which we have equivalent isotropy summands.
	\\
	\\The subsequent material is as follows. In Section \ref{preliminaries}, we provide some preliminaries, including some important formulas for our solutions to the PRP. In Section \ref{classification}, we provide a classification of all noncompact, simply connected $G/H$ with $G$ semi-simple, $H$ connected, and $G/H$ having two irreducible isotropy summands. In Section \ref{PRP}, we provide our solutions to the PRP with the exception of $SO_0(1,7)/G_2$ which contains isomorphic summands. Finally, in Section \ref{SO17}, we provide a complete solution for $SO_0(1,7)/G_2$.
	\\
	\\
	\textit{Acknowledgments: Special thanks is in order, first and foremost, to Michael Jablonski who introduced me to much of mathematics, including Riemannian Geometry and the Prescribed Ricci Curvature problem. Secondly, I would like to show my appreciation to Roi Docampo and Greg Muller who provided helpful insights in the use of Algebraic Geometry tools to solve equations involving polynomial and rational functions, tools that proved to be essential in the solving of the PRP in our exceptional case.}
	\\
	\\
	\textit{Author's Note: This work is a derivative of the author's Ph.D dissertation (\cite{mythesis}) and is primarily self-contained. There are, however, some of the more tedious of computations that are left out and cited to said work or to \cite{code}, in which some of the code developed for this work is housed and publicly available.}

	\section{2. Preliminaries}\label{preliminaries}
	\vspace{.2in}
	\noindent First, let us recall some basics of homogeneous geometry and Lie theory (See \cite{Besse}, \cite{Helg}, or \cite{doCarmo} for more details on the basics provided here). For a homogeneous space $G/H$, we have a \textit{reductive decomposition}, $\g = \h \oplus \p$ for which $\p$ is an $ad_\h$ invariant complement in $\g$. Moreover, the tangent space $T_{eH} G/H$ can be identified with $\p$ via:
	\begin{align*}
		x \to X_x = \frac{d}{dt}|_0 exp(tx).eH.
	\end{align*}
	\noindent From this identification, we can get a one-to-one correspondence between the $G$-invariant metrics on $G/H$ and the $Ad_H$ invariant inner products on $\p$ via:
	\begin{align*}
		g(X_x , Y_y) := (x,y).
	\end{align*}
	\noindent Moreover, if we further assume that $H$ is connected, we may equivalently think of $(.,.)$ as $ad_\h$ invariant. Now, from this identification of vector spaces and inner products, we are able to shift our perspective from calculating $ric$ at the geometric level to that of the Lie algebra level. This convenience we, like many others beforehand, will exploit thoroughly.
	\\
	\\In the pursuit of solving the PRP, let us recall a formula that will be useful in our consideration of the Ricci tensor. Restricting ourselves to the setting in which $G$ is semi-simple and $H$ is connected, from \cite{Besse} we have that 
	
	\begin{align}\label{ric}
		ric(x,y) =  - \frac{1}{2}\sum_{i}([x, e_i]_\p, [y, e_i]_\p) - \frac{1}{2}B(x,y) + \frac{1}{4}\sum_{i, j} ([e_i, e_j]_\p, x)([e_i , e_j]_\p , y)
	\end{align}

	\noindent where $\{e_i\}$ is an orthonormal basis with respect to our $ad_\h$ invariant inner product, $(.,.)$, on $\p$ in the reductive decomposition $\g = \h \oplus \p$. In \cite{Nik1}, working with such noncompact $G/H$ in which $H \subset K \subset G$ and $K$ is the maximal compact in semi-simple $G$, Nikonorov  provides a formula for the diagonal of the Ricci tensor that will prove to be useful for us. However, we first need to establish some notation, mostly following that of Nikonorov's work.
	\\
	\\
	Let $\p'$ and $\p''$ be such that $\g = \kk \oplus \p''$ is a Cartan decomposition and $\p'$ is an $ad_\h$ invariant complement in $\kk$ such that $\kk = \h \oplus \p'$. Now, let $\g = \h \oplus \p' \oplus \p''$ be decomposed into irreducible $ad_\h$ representations with $\p'' = \p_1 \oplus ... \oplus \p_u$ and $\p' = \p_{u+1} \oplus ... \oplus \p_v$ where $(.,.)$ and $\langle . , . \rangle = B(.,.)_{\p''} - B(.,.)_{\p'}$ are simultaneously diagonalized on $\p = \p' \oplus \p''$. Let the ordering $\p_1, ..., \p_u$ be given by $x_1 \leq ... \leq x_u$ where $(.,.)_{\p_i} = x_i \langle . , . \rangle_{\p_i}$. Similarly, we order the $\p_{u+1}, ..., \p_v$. Now, define $R_i$ to be $ric(.,.)_{\p_i} = R_i(.,.)_{\p_i}$. Let $b_i$ be 1 when $\p_i \subset \p'$ and -1 when $\p_i \subset \p''$, let $d_i$ be the dimension of $\p_i$, and let $\{e_i^\alpha\}$ be an orthonormal basis of $\p_i$ with respect to $\langle . , . \rangle$. By Lemma 2 in \cite{Nik1}, we have the following:
	
	\begin{align}\label{ricformula}
		R_i = \frac{b_i}{2x_i} + \frac{1}{4d_i}\sum_{1 \leq j, k \leq v} (\sum_{\alpha, \beta, \gamma} \langle [e_i^\alpha , e_j^\beta], e_k^\gamma \rangle^2) (\frac{x_i}{x_j x_k} - \frac{x_k}{x_i x_j} - \frac{x_j}{x_i x_k}).
	\end{align}
	
	\noindent It is important to observe that the $R_i$ describe the values of the $(1,1)$ Ricci tensor. Since we are primarily interested in working with the $(0,2)$ Ricci tensor, we will be interested in
	\begin{align}\label{riccidiagonals}
		r_i = x_i R_i.
	\end{align}	
	
	\begin{remark}
		Throughout this work, we will fix $\langle . , . \rangle = B_{\p''} - B_{\p'}$. This is a nice $ad_\h$ invariant inner product on $\p = \p' \oplus \p''$ for which $ad_x$ is skew-symmetric if $x \in \p'$ and symmetric if $x \in \p''$.
	\end{remark}
	
	\begin{remark}
		The formula for $r_i$ from Eqns. (\ref{ricformula}) and (\ref{riccidiagonals}) will be particularly useful as we examine the setting in which $G$ is simple. For $G$ semi-simple, thanks to our restricted setting, it will be enough for us to use the DeRham decomposition and Schur's lemma. See Section \ref{PRP} for more details.
	\end{remark}

	
	\section{3. Classification}\label{classification}
	
	
	\vspace{.2in}

	\noindent The following lemma is well known, but is included for completion as it is used in Theorem \ref{ClassificationTheorem} regarding the classification of two isotropy irreducible summand spaces.
	
	\begin{lemma}\label{IsoBetweenDuals}
		If $\g = \kk \oplus \p$ is a Cartan decomposition of a noncompact simple Lie algebra and $\g^\ast = \kk \oplus i\p$ is the dual to $\g$ (see Chapter V Section 2 and Chapterv VIII Section 5 of \cite{Helg}), then $\p$ and $i\p$ are isomorphic irreducible $ad_\kk$ representations.
	\end{lemma}
	\underline{Proof:}
	Consider the complexification $\g^\C = \kk \oplus \p \oplus i\kk \oplus i\p$ and restrict the scalars to be real. We may consider the real linear map $i: \p \to i\p$ and by definition of $\g^\C$, $[x, iv] = i[x,v]$, so we are done. \hfill $\blacksquare$
	\\
	\\Using the classification from \cite{Kerr} in the compact setting with two isotropy summands, we will classify spaces with two isotropy summands in the noncompact setting with the restriction that $G$ in $G/H$ is semi-simple, $G/H$ is simply connected, and $H$ connected. As with in the classification from \cite{Kerr}, the latter two restrictions are critical in our shift to the Lie algebra setting.
	\\
	\begin{theorem}\label{ClassificationTheorem}
		Let $G/H$ be simply connected with $G$ a connected semi-simple Lie group with no compact factors and $H \subset G$, a compact, connected subgroup. If $G/H$ has exactly two irreducible representations then $G/H$ is described by one of the following:
		\begin{enumerate}
			\item $G$ and $H$ have Lie algebras $\g = \g_1 \oplus \g_2$ and $\h = \kk_1 \oplus \kk_2$ where $\g_i$ is noncompact simple and $\kk_i$ is the maximal compact in $\g_i$. In this case, $G/H$ is a symmetric space. 
			\item $G$ has a simple Lie algebra $\g$, and $H$ has Lie algebra $\h \subsetneq \kk$ where $\kk$ is the maximal compact in $\g$. A classification of such $G/H$ is determined by the $(\g, \kk, \h)$ triple belonging to Tables \ref{TwoSummandsTable1} through \ref{TwoSummandsTable4}. In this case, $G/H$ is not a symmetric space.
		\end{enumerate}
	\end{theorem}

	\underline{Proof:}
	Let $\g$ be noncompact semi-simple and $\kk$ the maximal compact subalgebra of $\g$. If $\g$ is simple, then $G/K$ is an irreducible symmetric space, and the reductive decomposition of $\g$ is $\g = \kk \oplus \p$ with $\p$ being an irreducible $ad_\kk$ invariant complement to $\kk$ in $\g$. Therefore, if $\g$ is simple, we will have to consider an $\h \subsetneq \kk$, a case we will turn to after resolving the case in which $\g$ is semi-simple and not simple. 
	\\
	\\If $\g$ is not simple, then $G/K$ has a decomposition into irreducible factors coming from the DeRham decomposition. Thus, the reductive decomposition of $\g$ is $\g = \kk \oplus \p = \kk_1 \oplus ... \oplus \kk_n \oplus \p_1 \oplus ... \oplus \p_n$ in which $\p_i$ is the irreducible $ad_{\kk_i}$ invariant complement of $\kk_i$ in $\g_i$ and $[\g_i, \g_j] = 0$. Therefore, if $\g$ is semi-simple and not simple, we must restrict to $\g = \g_1 \oplus \g_2$ where $\g_1$ and $\g_2$ are simple ideals in $\g$. Moreover, we must have $\h = \kk_1 \oplus \kk_2$ where $\kk_1, \kk_2$ are the maximal compacts in $\g_1$ and $\g_2$, respectively. Indeed, if $\h \subsetneq \kk= \kk_1 \oplus \kk_2$, the maximal compact in $\g$, then the isotropy representation would have a (not necessarily irreducible) decomposition, $\p = \q \oplus \p_1 \oplus \p_2$ in which $\q$ is an $ad_\h$ invariant complement of $\h$ inside $\kk$. Thus, a complete description of spaces $G/H$ in which $G$ is semi-simple but not simple is given by $G_1/K_1 \times G_2/K_2$ in which $\g_i$ are noncompact simple and $\kk_i$ is the maximal compact inside $\g_i$. Such spaces are symmetric spaces. This gives us \textit{(1) }in the statement above.
	\\
	\\Now, let $\g$ be noncompact simple. As shown above, we must choose $\h \subsetneq \kk$ to obtain more than one irreducible summand. Since our reductive decomposition $\g = \h \oplus \p$ has a decomposition of $\p$ into irreducibles that is unique up to isomorphism, we may choose a decomposition of $\p$ that is convenient for us. Thus, we choose $\g = \h \oplus \p' \oplus \p''$ in which $\kk = \h \oplus \p'$ is a maximal compact in $\g$ and $\g = \kk \oplus \p''$ is a Cartan decomposition of $\g$. In general, $\p'$ and $\p''$ are not irreducible $ad_\h$ representations; however, since we are after those $G/H$ with two irreducible summands, $\p'$ and $\p''$ must be irreducible. One immediate consequence of our restricted setting is a restriction on the type of $\g$ we can have. Indeed, if $\g$ is the realification of a complex simple Lie algebra then $\g = \kk \oplus i\kk$ and $\g = \h \oplus \p' \oplus i\h \oplus i\p'$. Thus, our $\g$ must instead be a simple Lie algebra in which $\g^\C$ is complex simple.
	\\
	\\We now use the duality of symmetric spaces. We know that if we have a Cartan decomposition $\g = \kk \oplus \p''$ then $\g$ has a dual, $\g^\ast = \kk \oplus i\p''$ and $\g^\ast$ is compact with $G^\ast/K$ a compact irreducible symmetric space. Moreover, we know that $\p''$ and $i\p''$ are isomorphic $ad_\kk$ representations by Lemma \ref{IsoBetweenDuals}. Using this isomorphism of representations, we can pass from $\g = \h \oplus \p' \oplus \p''$ to $\g^\ast = \h \oplus \p' \oplus i\p''$ and vice versa. Therefore, for $\g$ simple, if $G/H$ is noncompact with two irreducible summands, we have $\g = \h \oplus \p' \oplus \p''$ and we will find a corresponding compact $G^\ast/H$ with two irreducible summands $\g^\ast = \h \oplus \p' \oplus i\p''$ with $\g^\ast$ compact simple. Using the complete list of compact $G^\ast/H$ with two irreducible summands in which $G^\ast$ is connected and simple given by Dickenson and Kerr in \cite{Kerr} (with corrections in Appendix A of \cite{IsoCorect} and Remark 6.1 of \cite{LauretTwoCorrection}), we can then get a complete list of noncompact $G/H$ with two irreducible summands. We now turn to how we can get that list.
	\\
	\\By observing the work in \cite{Kerr}, we note that the list in the noncompact setting will be smaller since in the compact setting one can have $G^\ast/H$ with two irreducible summands and $H$ maximal in $G^\ast$. In this case, $H$ will not be inside a $K \subset G^\ast$ such that $G^\ast/K$ is a compact irreducible symmetric space. Thus, such a $G^\ast/H$ could not be obtained from any noncompact $G/H$ with two irreducible summands using the duality. 
	\\
	\\In light of this, we must restrict ourselves in the compact setting to the case in which there is an intermediate subgroup $H \subset K \subset G^\ast$. Furthermore, since there are $H \subset K \subset G$ such that $K$ is maximal in $G^\ast$ but $G^\ast/K$ is isotropy irreducible and not symmetric (see tables 5 and 6 in Chapter 7 of \cite{Besse} for these), we must further restrict ourselves to those $H \subset K \subset G^\ast$ in which $G^\ast/K$ are irreducible symmetric spaces (which can be checked by the help of Chapter X of \cite{Helg}). We can then achieve our own list in the noncompact setting by the following procedure similar to the procedure used in Section 3 of \cite{ArrLaf}:
	\begin{enumerate}
		\item[\textbf{a.}] In the compact setting, select $\h \subset \kk \subset \g^\ast$ in which $(\g^\ast, \kk)$ is a pair associated with a compact irreducible symmetric space and $G^\ast/H$ has two irreducible summands. 
		\item[\textbf{b.}]  Use the duality to achieve $\h \subset \kk \subset \g$ such that $G/H$ has exactly two irreducible summands in the noncompact setting.
	\end{enumerate}
	A check of the lists in \cite{Kerr} (and corrections in \cite{IsoCorect} and \cite{LauretTwoCorrection}) shows that there are compact $G^\ast/H$ in which $H \subset K \subset G^\ast$ and $G^\ast/K$ is a compact irreducible symmetric space, but there are also some spaces in which $G^\ast/K$ is isotropy irreducible, but not symmetric. We wish to keep the former to use our duality, and ignore the latter. Thus, we can get our complete list in the noncompact setting by dualizing (in the sense of \textbf{b} above) each $G^\ast/H$ in the lists given in \cite{Kerr}, \cite{IsoCorect}, and \cite{LauretTwoCorrection} while ignoring the following items in the list: I.20, I.21, I.22, I.23, I.29, III.9, III.10, III.11, IV.3, IV.6, IV.13, IV.18, IV.30, IV.31, IV.32, IV.33, IV. 41, IV.42, IV.43, IV.44. \hfill $\blacksquare$
	
	\begin{remark}
		The tables below provide a correction to those found in \cite{mythesis}. An entire table was incidentally omitted in the final product of the prior work.
	\end{remark}

	\newpage
	
	\begin{table}[h!]
		\begin{center}
			\caption{}
			\label{tab:table}
			\begin{tabular}{l|c|c} \label{TwoSummandsTable1} 
				$(\g, \kk, \h)$ & Constraint & Label in \cite{Kerr}\\
				\hline
				\hline
				\\
				$(\so(m,k^2 +1), \so(m)\oplus\so(k^2 -1), \so(m)\oplus\su(k))$  & $k \geq 3$ & I.1\\
				$(\so(m, \frac{k}(k-1){2}), \so(m)\oplus\so(\frac{k(k-1)}{2}, \so(m)\oplus\so(k)))$  & $k \geq 7$ & I.2 \\
				$(\so(m, \frac{(k-1)(k+2)}{2}), \so(m)\oplus\so(\frac{(k-1)(k+2)}{2}), \so(m)\oplus\so(k))$  & $k\geq 5$ & I.3 \\
				$(\so(m, 2k^2 -k -1), \so(m)\oplus\so(2k^2 - k - 1) , \so(m)\oplus \sp(k))$  & $k \geq 3$ & I.4 \\
				$(\so(m, k(2k+1)), \so(m)\oplus\so(k(2k+1)), \so(m)\oplus\sp(k))$  & $k \geq 2$ & I.5 \\
				$(\so(m, 4k), \so(m)\oplus\so(4k), \so(m)\oplus\sp(1)\oplus\sp(k))$  & $k \geq 2$ & I.6 \\
				$(\so(m, 16), \so(m)\oplus\so(16), \so(m)\oplus\so(9))$ & NA & I.7 \\
				$(\so(m, 42), \so(m)\oplus\so(42), \so(m)\oplus \sp(4))$  & NA & I.8 \\
				$(\so(m,70), \so(m)\oplus\so(70), \so(m)\oplus\su(8))$ & NA & I.9 \\
				$(\so(m, 128), \so(m)\oplus\so(128), \so(m)\oplus\so(16))$  & NA & I.10 \\
				$(\so(m, 26), \so(m)\oplus\so(26), \so(m)\oplus\mathfrak{f}_4)$ & NA & I.11 \\
				$(\so(m,52), \so(m)\oplus\so(52), \so(m)\oplus\mathfrak{f}_4)$  & NA & I.12 \\
				$(\so(m, 78), \so(m)\oplus\so(78), \so(m)\oplus\mathfrak{e}_6)$  & NA & I.13 \\
				$(\so(m, 133), \so(m)\oplus\so(133), \so(m)\mathfrak{e}_7)$ & NA & I.14 \\
				$(\so(m,248), \so(m)\oplus\so(248), \so(m)\oplus \mathfrak{e}_8)$ & NA & I.15 \\
				$(\so(m,7), \so(m) \oplus \so(7), \so(m)\oplus\g_2)$  & $m = 1$, equal & I.16 \\
				$(\so(m,14), \so(m)\oplus\so(14), \so(m)\oplus\g_2)$  & NA & I.17 \\
				$(\so(m, 2k), \so(m)\oplus\so(2k), \so(m)\oplus\su(k)\oplus\R)$ & $k \geq 1$ & I.18 \\
				$(\so(3,5), \so(3)\oplus\so(5), \so(3)\oplus\so(3))$  & NA & I.19 \\
				$(\so(8,1), \so(8), \so(7))$  & NA & I.30 \\
				$(\so(m,8), \so(m)\oplus\so(8), \so(m)\oplus\so(7))$ & NA & \cite{IsoCorect} I.32 \\
				\hline
				\hline
				\\
				$(\so^\ast(2k), \su(k)\oplus\R , \su(k))$  & $k \geq 3$ & I.24 \\
				$(\so^\ast(2k), \su(k)\oplus\R, \so(k)\oplus\un(1))$  & $k\geq 3$ & I.25 \\
				$(\so^\ast(k^2 + k), \su(\frac{1}{2}(k^2 + k)\oplus\R), \su(k)\oplus\R)$  & $k\geq 3, \not= 4$ & I.26 \\
				$(\so^\ast(k^2 + k), \su(\frac{1}{2}(k^2 + k)\oplus\R), \su(k)\oplus\R))$  & $k \geq 5$ & I.27 \\
				$(\so^\ast(54), \su(27)\oplus\R , \mathfrak{e}_6\oplus\un(1))$ & NA & I.28 \\
				$(\so^\ast(32), \su(16)\oplus\R, \so(10)\oplus\un(1))$ & NA & \cite{IsoCorect} I.31 \\
				\hline
				\hline
			\end{tabular}
		\end{center}
	\end{table}
	\newpage

	\newpage
	
	\begin{landscape}
		
		\begin{table}[h!]
			\begin{center}
				\caption{}
				\label{tab:table}
				\begin{tabular}{l|c|c} 
					$(\g, \kk, \h)$ & Constraint & Label in \cite{Kerr}\\
					\hline
					\hline
					\\
					$(\su(\frac{n(n-1)}{2}, m), \su(\frac{n(n-1)}{2}) \oplus \su(m)\oplus \R, \su(n)\oplus\su(m) \oplus \R)$  & $n\geq 5$ & II.1 \\
					$(\su(\frac{n(n+1)}{2} , m), \su(\frac{n(n+1)}{2}) \oplus \su(m) \oplus \R, \su(n) \oplus \su(m) \oplus \R)$ & $n \geq 2$ & II.2 \\
					$(\su(27, m), \su(27)\oplus\su(m) \oplus \R, \mathfrak{e}_6 \oplus \su(m) \oplus \R)$  & NA & II.3 \\
					$(\su(16, m), \su(16) \oplus \su(m) \oplus \R, \so(10)\oplus\su(m)\oplus\R)$  & NA & II.4 \\
					$(\su(pq, m), \su(pq)\oplus\su(m)\oplus\R, \su(p)\oplus\su(q)\oplus\un(1)\oplus\su(m))$  & $p,q\geq 2$ & II.5 \\
					$(\su(n,m), \su(n)\oplus\su(m)\oplus\R, \so(n)\oplus\un(1)\oplus\su(m))$ & $n\geq 3$ & II.6 \\
					$(\su(n,m), \su(n)\oplus\su(m)\oplus\R, \su(n)\oplus\su(m))$  & $m,n\geq 2$ & II.7 \\
					$(\su(2n,m), \su(2n)\oplus\su(m)\oplus\R, \sp(n)\oplus\su(m)\oplus\R)$  & $n \geq 2$ & II.8 \\
					\hline
					\hline
					\\
					$(\su^\ast(14), \sp(7), \sp(3))$ & NA & II.9 \\
					$(\su^\ast(32), \sp(16), \so(12))$  & NA & II.10 \\
					$(\su^\ast(56), \sp(28), \mathfrak{e}_7)$  & NA & II.11 \\
					$(\su^\ast(4), \sp(2), \su(2))$ & $n\geq 5$ & II.1 \\
					\hline
					\hline
					\\
					$(\su^\ast(14), \sp(7), \sp(3))$ & NA & II.9 \\
					$(\su^\ast(32), \sp(16), \so(12))$  & NA & II.10 \\
					$(\su^\ast(56), \sp(28), \mathfrak{e}_7)$  & NA & II.11 \\
					$(\su^\ast(4), \sp(2), \su(2)))$  & NA & II.12 \\
					$(\su^\ast(20), \sp(10), \su(6))$  & NA & \cite{IsoCorect} II.15 \\
					\hline
					\hline 
					\\
					$(\sp(m,2), \sp(m)\oplus\sp(2), \sp(m)\oplus\su(2))$  & NA & III.1 \\
					$(\sp(m,7), \sp(m)\oplus\sp(7), \sp(m)\oplus\sp(3))$  & NA & III.2 \\
					$(\sp(m, 10), \sp(m)\oplus\sp(10), \sp(m)\oplus\su(6))$  & NA & III.3 \\
					$(\sp(m, 16), \sp(m)\oplus\sp(16), \sp(m)\oplus \so(12))$  & NA & III.4 \\
					$(\sp(m, 28), \sp(m)\oplus\sp(28), \sp(m)\oplus\mathfrak{e}_7)$ & NA & III.5 \\
					$(\sp(m,n) , \sp(m)\oplus\sp(n), \sp(m)\oplus\su(n)\oplus\R)$  & NA & III.6 \\
					$(\sp(m,n), \sp(m)\oplus\sp(n), \sp(m)\oplus\so(n)\sp(1))$  & $n \geq 3$ & III.7 \\
					\hline
					\hline
				\end{tabular}
			\end{center}
		\end{table}

	\end{landscape}

	\newpage

	\begin{table}[h!]
		\begin{center}
			\caption{}
			\label{tab:table}
			\begin{tabular}{l|c|c} 
				$(\g, \kk, \h)$  & Constraint & Label in \cite{Kerr}\\
				\hline
				\hline
				\\
				$(\sp(n, \R), \su(n)\oplus\R, \su(n))$  & $k \geq 3$ & III.8 \\
				$(\sp(2m, \R), \su(2m)\oplus\R , \sp(m)\oplus\un(1))$ & $m\geq 2$  & \cite{IsoCorect} III.12 \\
				\hline
				\hline
				\\
				$(\g_2(2),\so(4),\un(2))$ & $U(2)_3 \not\subset SU(3)$ & IV.1 \\
				$(\mathfrak{f}_4^{(4)}, \sp(3)\oplus\sp(1), \sp(3)\oplus\un(1))$ & NA & IV.2 \\
				$(\mathfrak{f}_4^{(-20)}, \so(9), \so(7)\oplus\so(2))$ & NA & IV.4 \\
				$(\mathfrak{f}_4^{(-20)}, \so(9), \so(6)\oplus\so(3))$  & NA & IV.5 \\
				\hline
				\hline
				\\
				$(\mathfrak{e}_6^{(-14)}, \so(10)\oplus\so(2), \so(10))$  & NA & IV.7 \\
				$(\mathfrak{e}_6^{(-14)}, \so(10)\oplus\so(2), \so(9)\oplus\so(2))$ & NA & IV.8 \\
				$(\mathfrak{e}_6^{(-14)}, \so(10)\oplus\so(2), \so(7)\oplus\so(3)\oplus\so(2))$ & NA & IV.9 \\
				$(\mathfrak{e}_6^{(-14)}, \so(10)\oplus\so(2), \so(5)\oplus\so(5)\oplus\so(2))$ & NA & IV.11 \\
				$(\mathfrak{e}_6^{(-14)}, \so(10)\oplus\so(2), \sp(2)\oplus\so(2))$ & NA & IV.12 \\
				$(\mathfrak{e}_6^{(2)}, \su(6) \oplus \su(2) , \su(6)\oplus\un(1))$  & NA & IV.14 \\
				$(\mathfrak{e}_6^{(2)}, \su(6)\oplus\su(2), \su(5)\oplus\R\oplus\su(2))$ & NA & IV.15 \\
				$(\mathfrak{e}_6^{(2)}, \su(6)\oplus\su(2), \so(6)\oplus\su(2))$ & NA & IV.16 \\
				$(\mathfrak{e}_6^{(2)}, \su(6)\oplus\su(2), \su(3)\oplus\su(2))$ & NA & IV.17 \\
				\hline
				\hline
				\\
				$(\mathfrak{e}_7^{(-25)}, \mathfrak{e}_6\oplus\so(2), \mathfrak{e}_6)$ & NA & IV.19 \\
				$(\mathfrak{e}_7^{(-25)}, \mathfrak{e}_6\oplus\so(2), \sp(4)\oplus\so(2))$ & NA & IV.20 \\
				$(\mathfrak{e}_7^{(-25)}, \mathfrak{e}_6\oplus\so(2), \g_2\oplus\so(2))$ & NA & IV.21 \\
				$(\mathfrak{e}_7^{(-25)}, \mathfrak{e}_6\oplus \so(2), 
				su(3)\oplus\so(2))$  & NA & IV.22 \\
				$(\mathfrak{e}_7^{(7)}, \su(8), \su(7)\oplus\R)$ & NA & IV.23 \\
				$(\mathfrak{e}_7^{(-5)}, \so(12)\oplus\sp(1), \so(12)\oplus\un(1))$ & NA & IV.24 \\
				$(\mathfrak{e}_7^{(-5)}, \so(12)\oplus\sp(1), \so(11)\oplus\sp(1))$  & NA & IV.25 \\
				$(\mathfrak{e}_7^{(-5)}, \so(12)\oplus\sp(1), \so(10)\oplus\so(2)\oplus\sp(1))$ & NA & IV.26 \\
				$(\mathfrak{e}_7^{(-5)}, \so(12)\oplus\sp(1), \so(9)\oplus\so(3)\oplus\sp(1))$  & NA & IV.27 \\
				$(\mathfrak{e}_7^{(-5)}, \so(12)\oplus\sp(1), \sp(7)\oplus\sp(5)\oplus\sp(1))$ & NA & IV. 28 \\
				$(\mathfrak{e}_7^{(-5)}, \so(12)\oplus\sp(1), \so(6)\oplus\so(6)\oplus\sp(1))$ & NA & IV.29 \\
				\hline
				\hline
				\\
			\end{tabular}
		\end{center}
	\end{table}

	\newpage

	\begin{table}[h!]
		\begin{center}
			\caption{}
			\label{tab:table}
			\begin{tabular}{l|c|c}\label{TwoSummandsTable4} 
				$(\g, \kk, \h)$ & Constraint & Label in \cite{Kerr}\\
				\hline
				\hline
				\\
				$(\mathfrak{e}_8^{(8)} , \so(16), \so(15))$  & NA & IV.34 \\
				$(\mathfrak{e}_8^{(8)} , \so(16), \so(15)\oplus\so(2))$  & NA & IV.35 \\
				$(\mathfrak{e}_8^{(8)} , \so(16), \so(13)\oplus\so(3))$ & NA & IV.36 \\
				$(\mathfrak{e}_8^{(8)} , \so(16), \so(11)\oplus\so(5))$  & NA & IV.38 \\
				$(\mathfrak{e}_8^{(8)} , \so(16), \so(10)\oplus\so(6))$  & NA & IV.39 \\
				$(\mathfrak{e}_8^{(8)} , \so(16), \so(9)\oplus\so(7))$  & NA & IV.40 \\
				$(\mathfrak{e}_8^{(8)} , \so(16), \so(9))$  & NA & \cite{LauretTwoCorrection} \\
				\hline
				\hline
				\\
				$(\mathfrak{e}_8^{(-24)}, \sp(1)\oplus\mathfrak{e}_7, \un(1)\oplus\mathfrak{e}_7)$ & NA & IV.45 \\
				$(\mathfrak{e}_8^{(-24)}, \sp(1)\oplus\mathfrak{e}_7, \sp(1)\oplus\su(8))$ & NA & IV.46 \\
				$(\mathfrak{e}_8^{(-24)}, \sp(1)\oplus\mathfrak{e}_7, \sp(1)\oplus\su(3))$  & NA & IV.47 \\
			\end{tabular}
		\end{center}
	\end{table}

	\begin{remark}\label{so17}
		As remarked at the end of Section 2 of \cite{Kerr} (and can be easily be checked with the corrections in \cite{IsoCorect} and \cite{LauretTwoCorrection}), the only $G^\ast/H$ with $\g^\ast = \h \oplus \p' \oplus i\p''$ in which $ \p' \simeq i\p''$ is $SO(8)/G_2$. Thus, in the noncompact setting, we similarly have that the only $(\g, \kk, \h)$ triple in which $\p' \simeq \p''$ is $(\so(1,7), \so(7), \g_2)$. The impact of this isomorphism on the Ricci curvature tensor turns out to be significant. Thus, the computation of the Ricci tensor and solving of the PRP in this setting is reserved to Section \ref{SO17}. However, utilizing the inequivalence of $\p'$ and $\p''$ for all other $(\g, \kk, \h)$ in the tables above, the solution to the PRP can be found with relatively little conceptual difficulty but a relatively large amount of tedium.
	\end{remark}

	\section{4. Prescribed Ricci Curvature}\label{PRP}
	
	\vspace{.2in}
	
	\noindent In the subsequent material, we are firmly settled into the setting in which $G/H$ has two irreducible summands and $G$ is semi-simple. For simple $\g$, then, we will be using the following notation for describing our decomposition of $\g$:

	\begin{align*}
		\g &= \h \oplus \p_1 \oplus \p_2  \text{ where }\\
		\kk &= \h \oplus \p_1  \text{ is a maximal compact in $\g$, and}\\
		\g &= \kk \oplus \p_2 \text{ is a Cartan decomposition of $\g$.} \\
	\end{align*}

	\noindent If, on the other hand, we have $\g$ semi-simple and not simple, then following Theorem \ref{ClassificationTheorem}, $\g = \g_1 \oplus \g_2$ with $\g_1$, $\g_2$ simple, and our decomposition of interest will come from the Cartan decomposition of each simple factor in the following manner:

	\begin{align*}
		\g &= \h \oplus \p_1 \oplus \p_2 \text{ where} \\
		\h &= \kk = \kk_1 \oplus \kk_2 \text{ is the maximal compact in $\g$, and} \\
		\g_1 &= \kk_1 \oplus \p_1 \text{ and } \g_2 = \kk_2 \oplus \p_2 \text{ are Cartan decompositions of each simple factor.} \\
	\end{align*}

	\noindent With this notation in place, we begin to investigate the kinds of $T$ we expect as solutions to $ric(.,.) = T(.,.)$ and $ric(.,.) = cT(.,.)$. Henceforth, unless there is a need for more specificity, we will reference $ric(.,.)$ as $ric$ and $T(.,.)$ as $T$.

	\begin{remark}
		We restrict ourselves to $c > 0$ since if $c < 0$, all we must do is consider the solutions to $ric(.,.) = cT(.,.)$ for $c> 0$ and negate $T(.,.)$.
	\end{remark}
	\begin{remark}
		As previously remarked, we fix $\langle .,. \rangle$ so that $\langle .,. \rangle = B_{\p''} - B_{\p'}$ where $B(.,.)$ is the Killing form for $\g = \h \oplus \p' \oplus \p''$. Note that for $G/H$ with two irreducible summands and $\g$ semi-simple and not simple, $\langle .,. \rangle = B_{\p''}$ as $\h = \kk$ in this case.
	\end{remark}
	
	\begin{lemma}\label{SchurTens}
		For $G/H$ with $G$ semi-simple noncompact and $G/H$ having two irreducible summands, every $(0,2)$ $ad_\h$ invariant tensor $T$ on the isotropy representation $\p = \p_1 \oplus \p_2$ is of the form $T = t_1\langle . , . \rangle_1 + t_2\langle . , . \rangle_2$ with the exception of $SO_0(1,7)/G_2$.
	\end{lemma} 
	\underline{Proof:}
	Since $T(.,.)$ is $ad_\h$ invariant, and since $\langle . , . \rangle$ is $ad_\h$ invariant, we know that in general, $T(x,y) = \langle \Phi x , y \rangle$ where $\Phi$ is symmetric and an $ad_\h$ equivariant map on $\p$. To prove the desired result, all we need to show is that $\Phi$ is necessarily diagonal. That is, we need to show that any $ad_\h$ intertwining map $\phi: \p_1 \to \p_2$ is 0.
	\\If $\g = \h \oplus \p_2 \oplus \p_1$ is simple, then, as noted in Remark \ref{so17}, since $\p_2 \not\simeq \p_1$ we have, by Schur's Lemma, that all $ad_\h$ intertwining maps $\p_1 \to \p_2$ are 0, providing the desired result. Similarly, when we have $\g$ semi-simple, by Theorem \ref{ClassificationTheorem}, we know $\g = \g_1 \oplus \g_2$ with $\h = \kk_1 \oplus \kk_2$ and $G/H = G_1/K_1\times G_2/K_2$. In this case, $ad_{\kk_i}$ acts irreducibly on $\p_i$ and trivially on $\p_j$ ($j \not= i$). This implies that no non-trivial, $ad_{\kk_1 \oplus \kk_2}$ intertwining maps $\p_2  \to \p_1$ exists, giving us the desired result. \hfill $\blacksquare$
	\\
	\\
	
	\begin{theorem}\label{ricequalsTsemisimple}
		If $\g$ is not simple then the only $T$ such that $ric = T$ is $T = ric_{\langle . , . \rangle}$ where $\langle . , . \rangle$ is our fixed inner product. Moreover, our $T < 0$ in this case. 
	\end{theorem}
	\underline{Proof:}
	Here, by the deRahm decomposition for symmetric spaces we have $ric = ric_1 + ric_2$ where $ric_1$ , $ric_2$ are the Ricci tensors for $G_1/K_1$ , $G_2/K_2$, respectively. Now, by $G_i/K_i$ being a noncompact irreducible symmetric space, we have that $ric_{\langle . , . \rangle_i} = \la_i \langle . , . \rangle_i$ for some $\la_i < 0$ for $i = 1,2$. Moreover, $ric_{c\langle . , . \rangle_i} = ric_{\langle . , . \rangle_i}$, and by Schur's Lemma, $\alpha\langle . , . \rangle_i$ exhausts all $ad_{\kk_i}$ inner products on $\p_i$ where $\alpha > 0$. Thus, for an arbitrary $ad_\h$ invariant inner product $(.,.)$, $ric_{(.,.)} = \la_1\langle . , . \rangle_1 + \la_2 \langle . , . \rangle_2 = ric_{\langle . , . \rangle}$. Therefore, $ric = T$ if and only if $T = ric_{\langle . , . \rangle}$, and so $T< 0$ by $\la_1, \la_2 < 0$. \hfill $\blacksquare$
	\\
	\\
	\begin{theorem}
		If $\g$ is not simple then $ric = cT$ for $c > 0$ if and only if $T$ is a scalar multiple of $ric_{\langle . , .\rangle}$ where $\langle . , . \rangle$ is our fixed inner product, and $c = \frac{T(x,x)}{ric_{\langle . , . \rangle}(x,x)}$ for some $x \in \p = \p_1 \oplus \p_2$. Moreover, our $T < 0$ in this case.
	\end{theorem}
	\underline{Proof:}
	The proof of Theorem \ref{ricequalsTsemisimple} allows us to easily see the solutions to $ric = cT$. Since $ric$ is given by $ric_{(.,.)} = ric_{\langle . , . \rangle} < 0$ for all $ad_\h$ invariant inner products, $(.,.)$, we have a solution to $ric = cT$  for $c > 0$ if and only if $T$ has $T < 0$ and is a scalar multiple of $ric_{\langle .,. \rangle}$ with $c = \frac{T(x,x)}{ric_{\langle . , . \rangle} (x,x)}$ for some $x \in \p = \p_1 \oplus \p_2$.
	\hfill $\blacksquare$
	\\
	\\Considering the PRP in the case where $\g$ is simple, the situation becomes more complex. One complexity is (as has been mentioned) that $\g = \so(1,7)$ must be handled differently in order to solve the PRP in its entirety (and will be handled in Section \ref{SO17}). Another complexity is how our formula for $ric$ can vary rather drastically with the bracket relation on $\p$. Due to these variations, we consider the solutions to $ric = T$ and $ric = cT$ in three different settings: first with $\p_1$ such that $[\p_1, \p_1] \not\subset \h$, second with $\p_1$ such that $[\p_1, \p_1] \subset \h$, and then third with $\p_1$ being a trivial $ad_\h$ representation. It turns out that the solutions in the second and third situations follow easily from the first and can, loosely, be thought of as specialized situations of the first one. Thus, we provide results for $ric = T$ and $ric = cT$ in the setting of $[\p_1, \p_1]\not\subset \h$, in Theorem \ref{ricequalsTgeneralp} and Theorem \ref{ricequalscTgeneral}, respectively, and then we provide corollaries describing the $[\p_1, \p_1] \subset \h$ and trivial representation settings. 
	\\
	\\Before providing solutions, though, we first provide formulas describing an arbitrary $(0,2)$ Ricci tensor in terms of an arbitrary $ad_\h$ inner product and Lie algebra data.
	
	\begin{remark}\label{Notation}
		For the subsequent statements and proofs, we utilize the notation from Eqn. \ref{ricformula}. Here, we have $d_1 = dim(\p_1)$, and $d_2 = dim(\p_2)$. Moreover, we set $p_1 = \sum_{\alpha,\beta,\gamma} \langle [e_1^\alpha, e_1^\beta], e_1^\gamma \rangle^2$ and $p_2 = \sum_{\alpha,\beta,\gamma} \langle [e_2^\alpha, e_1^\beta], e_2^\gamma \rangle^2$. Recall that $\{e_1^\alpha\}$ and $\{e_2^\alpha\}$ are understood to be orthonormal bases with respect to $\langle . , . \rangle$, our fixed inner product on $\p_1$ and $\p_2$, respectively.
	\end{remark}
	\begin{lemma}\label{ricfrommetric}
		Let $G/H$ being a noncompact space with two irreducible summands in which $\g$ is simple with $\g = \h \oplus \p_1 \oplus \p_2$. In such a setting we have the following formulas for $r_1, r_2$ defining $ric(.,.)= r_1\langle . , . \rangle_1 + r_2\langle . , . \rangle_2$:
		\begin{align*}
			r_1 &= \dfrac{2d_1 - p_1 - 2p_2}{4d_1} + \frac{p_2}{4d_1}\left(\frac{x_1}{x_2}\right)^2 > 0 \\
			r_2 &= \frac{-1}{2} - \frac{p_2}{2d_2}\frac{x_1}{x_2} < 0 .\\
		\end{align*}
	\end{lemma}
	\underline{Proof:}
	Following the expressions provided by Nikonorov in \cite{Nik1} (as discussed in Section \ref{preliminaries}), we define an arbitrary $ad_\h$ invariant inner product $(.,.) = x_1\langle.,.\rangle_1 + x_2\langle.,.\rangle_2$ where $x_i > 0$ and the Ricci tensor for that inner product, $ric = R_1(.,.)_1 + R_2(.,.)_2$ where $R_1, R_2 \in \R$. From Eqn.\ref{ricformula} we have:
	$$R_1 = \frac{-1}{2x_1} + \frac{1}{4d_1}\sum_{1 \leq j,k \leq 2}(\sum_{\alpha,\beta,\gamma} \langle [e_1^\alpha, e_j^\beta], e_k^\gamma \rangle^2)(\frac{x_1}{x_j x_k} - \frac{x_k}{x_1 x_j} - \frac{x_j}{x_1 x_k})$$
	$$R_2 = \frac{1}{2x_2} + \frac{1}{4d_2}\sum_{1 \leq j,k \leq 2}(\sum_{\alpha,\beta,\gamma} \langle [e_2^\alpha, e_j^\beta], e_k^\gamma \rangle^2)(\frac{x_2}{x_j x_k} - \frac{x_k}{x_2 x_j} - \frac{x_j}{x_2 x_k}).$$
	\\
	\\The first step will be to generate simplified formulas for $R_1$ and $R_2$ in our setting dependent upon $x_1$ and $x_2$ (along with terms $p_1$ and $p_2$ dependent upon the bracket of $\g$ on $\p$). From there, we will determine our $r_1$ and $r_2$.
	\\
	\\By the Cartan decomposition properties and the (skew) symmetry of $ad_{e_i}$ with respect to $\langle . , . \rangle $, we get that in $R_1$, our term $\langle [e_1^\alpha, e_j^\beta], e_k^\gamma \rangle^2$ has:
	\begin{align*}
		\langle [e_1^\alpha, e_1^\beta], e_2^\gamma \rangle^2 = \langle [e_1^\alpha, e_2^\beta], e_1^\gamma \rangle^2 &=0 .
	\end{align*}
	
	Similarly, in $R_2$ our $\langle [e_2^\alpha, e_j^\beta], e_k^\gamma \rangle^2$ term has:
	\begin{align*}
		\langle [e_2^\alpha, e_1^\beta], e_2^\gamma \rangle^2 &= \langle [e_2^\alpha, e_2^\beta], e_1^\gamma \rangle^2 \\
		\langle [e_2^\alpha , e_2^\beta], e_2^\gamma \rangle^2 &= \langle [e_2^\alpha, e_1^\beta], e_1^\gamma \rangle^2=0. \\
	\end{align*}
	Letting $p_1 = \sum_{\alpha,\beta,\gamma} \langle [e_1^\alpha, e_1^\beta], e_1^\gamma \rangle^2$, and $p_2 = \sum_{\alpha,\beta,\gamma}\langle [e_2^\alpha, e_2^\beta], e_1^\gamma \rangle^2 = \sum_{\alpha,\beta,\gamma} \langle [e_2^\alpha, e_1^\beta], e_2^\gamma \rangle^2 = \sum_{\alpha,\beta,\gamma} \langle [e_1^\alpha, e_2^\beta], e_2^\gamma \rangle^2$, we get the following formulas for $R_1$ and $R_2$:
	$$R_1 = \frac{1}{2x_1} + \frac{1}{4d_1}\left[p_2\left(\frac{x_1}{x_2 x_2} - \frac{x_2}{x_1 x_2} - \frac{x_2}{x_1 x_2} \right) + p_1\left(\frac{x_1}{x_2 x_1} - \frac{x_1}{x_1 x_1} - \frac{x_1}{x_1 x_1} \right) \right].$$
	$$R_2 = \frac{-1}{2x_2} + \frac{p_2}{4d_2}\left[\left(\frac{x_2}{x_2 x_1} - \frac{x_1}{x_2 x_2} - \frac{x_2}{x_2 x_1}\right) + \left(\frac{x_2}{x_1 x_2} - \frac{x_2}{x_2 x_1} - \frac{x_1}{x_2 x_2}\right)\right]$$
	Simplifying both terms and getting a common denominator we get the following:
	\begin{align*}
		R_1 &= \dfrac{p_2(x_1^2 - 2x_2) + (2d_1 - p_1)x_2^2}{4d_1x_1x_2^2} \\
		\\
		&= \dfrac{(2d_1 - p_1 - 2p_2)x_2^2 + p_2x_1^2}{4d_1x_1x_2^2}. \\
		\\
		R_2 &= \dfrac{-d_2x_2 - p_2 x_1}{2d_2 x_2^2} \\
	\end{align*}
	By Lemma 1 in \cite{Nik1}, we have $p_1 + p_2 \leq d_1$ and $2p_2 \leq d_2$, with equality only when $\p_1$ and $\p_2$, respectively, are trivial representations for $ad_\h$. Moreover, since $p_2 = \sum_{\alpha,\beta,\gamma}\langle [e_2^\alpha, e_2^\beta], e_1^\gamma \rangle^2$, we know that $p_2 > 0$ since by the Cartan decomposition properties $[\p_2, \p_2] = \kk \supset \p_1$.  We can thus observe that $R_1 > 0$ and $R_2 < 0$.
	\\
	\\We now relate $ric$ to the background inner product, placing ourselves in the $(0,2)$ tensor setting. Since $(.,.) = x_1\langle . , . \rangle_1 + x_2\langle . , . \rangle_2$ and $ric(.,.) = R_1(.,.)_1 + R_2(.,.)_2$, we can write $ric(.,.) = x_1R_1 \langle . , . \rangle_1 + x_2R_2 \langle .,. \rangle$. With $r_1 = x_1R_1$ and $r_2 = x_2R_2$ and we have the following:
	\begin{align}\label{R1R2}
		r_1 &= x_1\dfrac{(2d_1 - p_1 - 2p_2)x_2^2 + p_2x_1^2}{4d_1x_1x_2^2} \notag\\
		&= \dfrac{2d_1 - p_1 - 2p_2}{4d_1} + \frac{p_2}{4d_1}\left(\frac{x_1}{x_2}\right)^2 \notag \\
		\notag\\
		r_2 &= x_2\dfrac{-d_2x_2 - p_2 x_1}{2d_2 x_2^2} \notag\\
		&=  \frac{-1}{2} - \frac{p_2}{2d_2}\frac{x_1}{x_2}
	\end{align} 
	Note that $r_1 = x_1R_1 >0$ since $x_1 > 0$ and $R_1 > 0$; likewise, $r_2 = x_2 R_2 < 0$ since $x_2 > 0$ and $R_2 < 0$. Thus, we have our desired result.\hfill $\blacksquare$
	\\
	\\
	\begin{lemma}\label{pvalueinequalities}
		Let $p_1, p_2, d_1 ,d_2$, and $\g = \h \oplus \p_2 \oplus \p_1$ be defined as in Remark \ref{Notation} and Lemma \ref{ricfrommetric}. In general, we have $$\dfrac{2d_1-p_1-2p_2}{4d_1} \geq 0.$$
		If we assume that $[\p_1, \p_1] \not\subset \h$, then we have $$\dfrac{2d_1-p_1-2p_2}{4d_1} > 0.$$
		If we assume that $[\p_1, \p_1] \subset \h$, then we have $$\dfrac{2d_1-p_1-2p_2}{4d_1} = \dfrac{d_1-p_2}{2d_1} \geq 0$$ with equality if and only if $\p_1$ is a trivial $ad_\h$ representation.
	\end{lemma}
	\underline{Proof:}
	As was mentioned in the proof of Lemma \ref{ricfrommetric}, by Lemma 1 in \cite{Nik1}, $p_1 + p_2 \leq d_1$ with equality if and only if $\p_1$ is a trivial representation. By $p_1 + p_2 \leq d_1$ we are able to conclude that $\dfrac{2d_1-p_1-2p_2}{4d_1} \geq 0$ for general $\p_1$, proving the first claim.
	\\
	\\If $[\p_1, \p_1] \not \subset \h$ then we know that $p_1 =\sum_{\alpha,\beta,\gamma} \langle [e_1^\alpha, e_1^\beta], e_1^\gamma \rangle^2 > 0$. Thus, we have that $$2d_1 - p_1 - 2p_2 > 2d_1 - 2p_1 - 2p_2 \geq 0$$ and we are able to conclude that $\dfrac{2d_1-p_1-2p_2}{4d_1} > 0$, proving the second claim.
	\\
	\\If $[\p_1, \p_1] \subset \h$ then we know that $p_1 = 0$ and we get $\dfrac{2d_1-p_1-2p_2}{4d_1} = \dfrac{d_1-p_2}{2d_1} \geq 0$ with $p_2 = d_1$ if and only if $\p_1$ is trivial, as mentioned above. \hfill $\blacksquare$
	\\
	\\
	\begin{theorem}\label{ricequalsTgeneralp}
		Let $\g = \h \oplus \p_1 \oplus \p_2$ be noncompact simple with $[\p_1, \p_1] \not\subset \h$ (except $\g = \so(1,7)$). For $G/H$ in this case, $ric = T$ has a solution $T$ if and only if $$t_1 = \frac{d_2^2}{d_1p_2}t_2^2 + \frac{d_2^2}{d_1p_2}t_2 + \dfrac{p_2(2d_1 - p_1 - 2p_2) + d_2^2}{4d_1p_2} \text{ with } t_2 \in (-\infty, \frac{-1}{2}).$$
	\end{theorem}
	\underline{Proof:}
	The goal of this proof is to determine sufficient and necessary conditions on $T(.,.) = t_1\langle . , . \rangle_1 + t_2\langle . , . \rangle_2$ such that $ric = T$ for some $ad_\h$ invariant inner product. Using $r_1$ and $r_2$ defining $ric(.,.)$ in Lemma \ref{ricfrommetric}, we find our solution utilizing the following system of equations:
	$$\begin{cases}
		r_1 &= t_1 \\
		r_2 &= t_2.
	\end{cases}$$
	Note that by $r_1 > 0$ and $r_2 < 0$, we know that $t_1 > 0$ and $t_2 < 0$. Now, plugging into $r_1$ and $r_2$ we have the following:
	\begin{align*}
		\dfrac{2d_1 - p_1 - 2p_2}{4d_1} + \frac{p_2}{4d_1}\left(\frac{x_1}{x_2}\right)^2 &= t_1 \\
		\frac{-1}{2} - \frac{p_2}{2d_2}\frac{x_1}{x_2} &= t_2.
	\end{align*}
	Let $\la = \frac{x_1}{x_2}$ and observe that $\la$ can take on any positive value, implying that $t_2$ can take on any value in $(-\infty, \frac{-1}{2})$. Our approach is as follows. We use the equation with $t_2$ to solve for $\la$ and then we substitute $\la$ into the equation with $t_1$, providing an equation of $t_1$ in terms of $t_2$. Once we have that, we will know that for any $t_2 \in (-\infty, \frac{-1}{2})$, we can get a $t_1$ such that $ric = T$ has a solution, providing sufficient and necessary conditions as desired.
	\\
	\\Thus, we have:
	\begin{align*}
		\frac{-1}{2} - \frac{p_2}{2d_1}\la &= t_2 \\
		\la &= \frac{-2d_2}{p_2}(t_2 + \frac{1}{2}) \\
		\\
	\end{align*}
	\begin{align*}
		\dfrac{2d_1 - p_1 - 2p_2}{4d_1} + \frac{p_2}{4d_1}\left(\la\right)^2 &= t_1\\
		\dfrac{2d_1 - p_1 - 2p_2}{4d_1} + \frac{p_2}{4d_1}\left(\frac{-2d_2}{p_2}(t_2 + \frac{1}{2})\right)^2 &= t_1 \\
		\dfrac{2d_1 - p_1 - 2p_2}{4d_1} + \frac{p_2}{4d_1}\left(\frac{4d_2^2}{p_2^2}(t_2^2 + t_2 + \frac{1}{4})\right) &= t_1 \\ 
		\dfrac{2d_1 - p_1 - 2p_2}{4d_1} + \frac{d_2^2}{d_1p_2}\left(t_2^2 + t_2 + \frac{1}{4}\right) &= t_1 \\
		\frac{d_2^2}{d_1p_2}t_2^2 + \frac{d_2^2}{d_1p_2}t_2 + \dfrac{p_2(2d_1 - p_1 - 2p_2) + d_2^2}{4d_1p_2} &= t_1
	\end{align*}
	Thus, the solutions to $ric = T$ are given by $$t_1 = \frac{d_2^2}{d_1p_2}t_2^2 + \frac{d_2^2}{d_1p_2}t_2 + \dfrac{p_2(2d_1 - p_1 - 2p_2) + d_2^2}{4d_1p_2} \text{ with } t_2 \in (-\infty, \frac{-1}{2}).$$
	\hfill $\blacksquare$
	
	\begin{remark}
		In addition to knowing the $T$ for which we have solutions to $ric = T$, in the course of the proof, the metric is also made clear. By the equality $\la = \frac{-2d_2}{p_2}(t_2 + \frac{1}{2})$, we know that for a given $G/H$ and $T = t_1\langle . , . \rangle + t_2 \langle .,. \rangle $, the $(.,.) = x_1 \langle . , . \rangle_1 + x_2\langle . , . \rangle_2$ for which $ric_{(.,.)} = T$ is one in which $\frac{x_1}{x_2} = \frac{-2d_2}{p_2}(t_2 + \frac{1}{2})$.
		\\
	\end{remark}
	
	\begin{theorem}\label{ricequalscTgeneral}
		Let $\g = \h \oplus \p_1 \oplus \p_2$ be noncompact simple with $[\p_1, \p_1] \not\subset \h$ (except $\g = \so(1,7)$). For $G/H$ in this case, the equation $ric = cT$ for $c> 0$ has a solution if and only if $(t_1, t_2)$ is a pair satisfying $$\frac{t_1}{t_2} \leq \frac{d_2^2}{d_1p_2} - \frac{d_2}{d_1}\sqrt{\frac{d_2^2}{p_2^2} + \dfrac{2d_1 - p_1 - 2p_2}{p_2}} \leq 0$$
		where $t_1 < 0$ and $t_2 > 0$. (Equality on the right only occurs when $\p_1$ is trivial. See Corollary \ref{ricequalscTtrivialp2} for that setting).
		\\
		\\When the above inequality is satisfied, there is always one solution which can be obtained by $c_+$ (where $c_+$ takes the $+$ in (\ref{cval}) below) and $(x_1, x_2)$, the pair unique up to scaling given by $\frac{x_2}{x_1} = \frac{-2d_2}{p_2}(c_+ t_2 + \frac{1}{2}).$ 
		\\
		\\In addition to this one solution, there is a second solution if and only if our pair $(t_1, t_2)$ with $t_1 < 0$ and $t_2 > 0$ satisfies
		$$-\dfrac{2d_1 - p_1 - 2p_2}{2d_1} < \frac{t_1}{t_2} < \frac{d_2^2}{d_1p_2} - \frac{d_2}{d_1}\sqrt{\frac{d_2^2}{p_2^2} + \dfrac{2d_1 - p_1 - 2p_2}{p_2}}.$$
		The second solution is obtained by $c_-$ (where $c_-$ takes the $-$ in \ref{cval} below) and $(x_1, x_2)$, the pair unique up to scaling given by $\frac{x_1}{x_2} = \frac{-2d_2}{p_2}(c_- t_2 + \frac{1}{2}).$
		\begin{align}\label{cval}
			c = \dfrac{-\left(\frac{d_2^2}{d_1p_2}t_2 - t_1\right) \pm \sqrt{\left(\frac{d_2^2}{d_1p_2}t_2 - t_1\right)^2 - 4\left(\frac{d_2^2}{d_1p_2}t_2^2\right)\left(\frac{d_2^2}{4p_2d_1} + \dfrac{2d_1-p_1-2p_2}{4d_1}\right)}}{2\frac{d_2^2}{d_1p_2}t_2^2} 
		\end{align}
		Moreover, $c_+ = c_-$ when the discriminant in $c$ is zero, and this happens precisely when $$\frac{t_1}{t_2} = \frac{d_2^2}{d_1p_2} - \frac{d_2}{d_1}\sqrt{\frac{d_2^2}{p_2^2} + \dfrac{2d_1 - p_1 - 2p_2}{p_2}}.$$
	\end{theorem}
	
	\begin{remark}
		If $\g = \so(1,7)$ then Theorems \ref{ricequalsTgeneralp} and \ref{ricequalscTgeneral} still hold so long as we restrict ourselves to those inner products in which $(\p_1, \p_2) = 0$.
	\end{remark}
	
	\includegraphics[width=\textwidth]{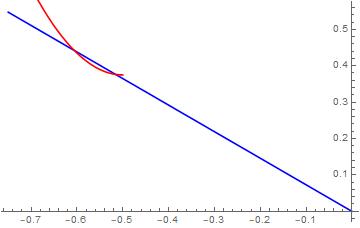}
	\noindent The above image was produced in Mathematica (\cite{Mathematica}) for the image of $ric$ (in red) with $G/H = SO_0(1,7)/G_2$ when we restrict our $ad_{\g_2}$ inner products to those in which $(\p_2, \p_1) = 0$. In this case $d_1 = d_2 = 7$ and $p_1 = p_2 = \frac{7}{6}$. Since $ric = cT$ is looking for $T$ produced by taking the image of $ric$ and multiplying by some constant, all such $T$ can be understood by what is called the \textit{cone} of the image of $ric$ (See Exercise 16 in Chapter 8 of Section 2 in \cite{Cox}). The cone of the image of $ric$ is the collection of lines through the origin that intersect the image of $ric$. Since we concern ourselves with $c> 0$ specifically, we ignore one half of the line through the origin depending on the point on the image of $ric$ under consideration. In our case $t_1 > 0$ and $t_2 < 0$, so we would always ignore the half of the line at the origin and in the fourth quadrant if we orient our plane with $t_2$ being the horizontal axis and $t_1$ the vertical. The line provided above helps to illustrate a subset of points describing $T$ in $ric = cT$. Additionally, the line illustrates the need for more than one $c$ value depending on how the line intersects $T$ in certain cases.
	
	\begin{remark}\label{specialpvalues}
		A remark is warranted before we begin our proof. In the following, we will be using Lemma \ref{pvalueinequalities} to come to conclusions regarding the existence of solutions to $ric = cT$. For the purposes of proving the corollaries to follow in a simpler fashion, we use the general $\p_1$ setting of the lemma unless forced to assume $[\p_1, \p_1] \not\subset \h$. We will see that the only setting that does more than a minor formula change for the $c$ and the $(t_1, t_2)$ is the setting in which $\p_1$ is a trivial representation. As will be shown, this is the case because $\p_1$ trivial is the only setting in which $\dfrac{2d_2 - p_2 - 2p_1}{p_1} = 0$. When the need for strict inequality in Lemma \ref{pvalueinequalities} occurs, we are sure to make note of it and point the reader to Corollary \ref{ricequalscTtrivialp2}.
	\end{remark}

	\underline{Proof:}
	The goal of this proof is to find necessary and sufficient conditions on $T = t_1 \langle . , . \rangle_1 + t_2\langle . , . \rangle_2$ such there there are solutions to $ric = cT$. In the course of the proof, we also provide what $c > 0$ and $ad_\h$ inner product (determined by the pair $(x_1, x_2)$) to expect for a given $(t_1, t_2)$ for which there is a solution. As before, using $r_1$ and $r_2$ defining $ric(.,.)$ in Lemma \ref{ricfrommetric}, we find our solution utilizing the following system of equations:
	$$\begin{cases}
		r_1 &= ct_1 \\
		r_2 &= ct_2.
	\end{cases}$$
	Plugging in for $r_1, r_2$, we get the following equations:
	\begin{align*}
		\dfrac{2d_1 - p_1 - 2p_2}{4d_1} + \frac{p_2}{4d_1}\left(\frac{x_1}{x_2}\right)^2 &= ct_1 \\
		\frac{-1}{2} - \frac{p_2}{2d_2}\frac{x_1}{x_2}  &= ct_2.
	\end{align*}
	Note that once again, by $r_1, c > 0$ and $r_2 < 0$ we have $t_1 > 0$ and $t_2 < 0$. Our approach is similar to that of Theorem \ref{ricequalsTgeneralp}, but since we need more information than just what $t_1, t_2$ satisfy the equation, there are some differences that make the solutions more tedious to find.
	\\
	\\In this case we are looking for conditions on the pair $(t_1, t_2)$ that are sufficient and necessary to the existence of a $c> 0$ and $(x_1, x_2)$ such that we have a solution to the given system of equations. To do so, we again set $\la = \frac{x_1}{x_2}$, and solve for $\la$ in terms of $c$ and $t_2$. We then get an equation for $c$ in terms of $t_1$ and $t_2$, and using the condition that $c> 0$, we obtain sufficient and necessary conditions on $(t_1, t_2)$ such that we get a $c>0$. Then, we determine which subset of those $(t_1, t_2)$ giving us $c > 0$ satisfy $\la > 0$. This will provide us with solutions to the given system of equations, providing us with the desired result.
	\begin{align*}
		\frac{-1}{2} - \frac{p_2}{2d_2}\frac{x_1}{x_2} &= ct_2  \\
		\frac{-1}{2} - \frac{p_2}{2d_2}\la &= ct_2 \\
		\la&= \frac{-2d_2}{p_2}\left(ct_2 + \frac{1}{2}\right) \tag{$\star$} \label{c2.2}\\
	\end{align*}
	\begin{align*}
		\dfrac{2d_1 - p_1 - 2p_2}{4d_1} + \frac{p_2}{4d_1}\left(\frac{x_1}{x_2}\right)^2 &= ct_1 \\
		\dfrac{2d_1 - p_1 - 2p_2}{4d_1} + \frac{p_2}{4d_1}\la^2 &= ct_1 \\
		\dfrac{2d_1 - p_1 - 2p_2}{4d_1} + \frac{p_2}{4d_1}\left(\frac{-2d_2}{p_2}\left(ct_2 + \frac{1}{2}\right)\right)^2 &= ct_1 \\
		\dfrac{2d_1 - p_1 - 2p_2}{4d_1} + \frac{p_2}{4d_1}\left(\frac{4d_2^2}{p_2^2}\left(c^2t_2^2 + ct_2 + \frac{1}{4}\right)\right) &= ct_1  \\
		\left(\frac{d_2^2}{d_1p_2}t_2^2\right)c^2 + \left(\frac{d_2^2}{d_1p_2}t_2 - t_1\right)c + \frac{d_2^2}{4p_2d_1} + \dfrac{2d_1 - p_1 - 2p_2}{4d_1} &= 0 \tag{$\star \star$}\\
	\end{align*}
	Observing that the above equation is quadratic in $c$, we solve for $c$ using the quadratic formula.
	\\
	\begin{align}
		c = \dfrac{-\left(\frac{d_2^2}{d_1p_2}t_2 - t_1\right) \pm \sqrt{\left(\frac{d_2^2}{d_1p_2}t_2 - t_1\right)^2 - 4\left(\frac{d_2^2}{d_1p_2}t_2^2\right)\left(\frac{d_2^2}{4p_2d_1} + \dfrac{2d_1-p_1-2p_2}{4d_1}\right)}}{2\frac{d_2^2}{d_1p_2}t_2^2}  \label{cvals} 
		\\\notag
	\end{align}
	So long as $c > 0$ and the resulting $\la > 0$, any $t_1, t_2, c$ satisfying ($\star \star$) above provides a solution to $ric = cT$. Thus, we use $c$ above to find what conditions on $t_1$ and $t_2$ are necessary and sufficient for $c > 0$ and then we use ($\star$) to determine what conditions are necessary and sufficient for $\la > 0$. Once we have those, we will have all the solutions to ($\star \star$) and thus all the solutions to $ric = cT$. There are multiple steps here with (in the end) more than one possible solution in certain settings. For this reason, we finish the proof with a set of claims:
	\begin{itemize}
		\item[\textbf{Claim 1:}] $c > 0$ if and only if $c$ is real
		\item[\textbf{Claim 2:}] $c$ is real if and only if $(t_1, t_2)$ satisfies (\ref{tvalsclaim})
		\begin{align}
			\frac{t_1}{t_2} \leq \frac{d_2^2}{d_1p_2} - \frac{d_2}{d_1}\sqrt{\frac{d_2^2}{p_2^2} + \dfrac{2d_1 - p_1 - 2p_2}{p_2}} \label{tvalsclaim}
		\end{align}
		\item[\textbf{Claim 3:}] For any $(t_1, t_2)$ satisfying (\ref{tvalsclaim}), there is one solution which can be obtained by $c_+$ and $\la = \frac{-2d_2}{p_2}\left(c_+t_2 + \frac{1}{2}\right)$ where $c_+$ takes the $+$ in (\ref{cvals})
		\item[\textbf{Claim 4:}] For any $(t_1, t_2)$ satisfying (\ref{tvalssecondsolution}), there is a second solution which can be obtained by $c_-$ and $\la = \frac{-2d_2}{p_2}\left(c_-t_2 + \frac{1}{2}\right)$ where $c_-$ takes the $-$ in (\ref{cvals}). 
		\\
		\begin{align}\label{tvalssecondsolution}
			-\dfrac{2d_1 - p_1 - 2p_2}{2d_1} < \frac{t_1}{t_2} < \frac{d_2^2}{d_1p_2} - \frac{d_2}{d_1}\sqrt{\frac{d_2^2}{p_2^2} + \dfrac{2d_1 - p_1 - 2p_2}{p_2}}
			\\\notag
		\end{align}
	\end{itemize}
	\noindent\textbf{\underline{Claim 1:}} $c > 0$ if and only if $c$ is real
	\\
	\\\underline{Proof of \textbf{Claim 1}:}
	First recall that $t_1 > 0$, $t_2 < 0$, and $\dfrac{2d_1-p_1-2p_2}{4d_1} \geq 0$ by the general setting in Lemma \ref{pvalueinequalities}. These inequalities allows us to also conclude that $$4\left(\frac{d_2^2}{d_1p_2}t_2^2\right)\left(\frac{d_2^2}{4p_2d_1} + \dfrac{2d_1-p_1-2p_2}{4d_1}\right) > 0.$$ 
	\\
	\\These four inequalities allow us to see that $c > 0$ if and only if $c$ is real as the numerator and denominator in $c$ will always be positive. Indeed, the denominator is always positive since every term is positive. The numerator being positive requires some more detailed checking. First, recall the numerator of $c$:
	\begin{align*}
		-\left(\frac{d_2^2}{d_1p_2}t_2 - t_1\right) \pm \sqrt{\left(\frac{d_2^2}{d_1p_2}t_2 - t_1\right)^2 - 4\left(\frac{d_2^2}{d_1p_2}t_2^2\right)\left(\frac{d_2^2}{4p_2d_1} + \dfrac{2d_1-p_1-2p_2}{4d_1}\right)}
	\end{align*}
	which is always positive if we take the plus sign and $c$ is real since $t_1 > 0$ and $t_2 < 0$. If we take the minus sign, then we need the term on the right to have smaller magnitude than the term on the left, which is true since the inequality $4\left(\frac{d_2^2}{d_1p_2}t_2^2\right)\left(\frac{d_2^2}{4p_2d_1} + \dfrac{2d_1-p_1-2p_2}{4d_1}\right) > 0$ implies that

	$$\sqrt{\left(\frac{d_2^2}{d_1p_2}t_2 - t_1\right)^2 - 4\left(\frac{d_2^2}{d_1p_2}t_2^2\right)\left(\frac{d_2^2}{4p_2d_1} + \dfrac{2d_1-p_1-2p_2}{4d_1}\right)} < \left|\left(\frac{d_2^2}{d_1p_2}t_2 - t_1\right)\right|.$$
	
	\vspace{.2in}
	
	\noindent Therefore, we have a $c > 0$ just by finding when $c$ is real. This concludes the proof of \textbf{Claim 1}.
	\\
	\\We now focus on finding conditions on $(t_1,t_2)$ when $c$ is real. 
	\\
	\\\textbf{\underline{Claim 2:}} $c$ is real if and only if $(t_1, t_2)$ satisfies (\ref{tvalsclaim})
	\\
	\\\underline{Proof of \textbf{Claim 2}:}
	To prove this claim, we look for $t_1, t_2$ such that the discriminant in $c$ is non-negative:
	\begin{align*}
		\left(\frac{d_2^2}{d_1p_2}t_2 - t_1\right)^2 - 4\left(\frac{d_2^2}{d_1p_2}t_2^2\right)\left(\frac{d_2^2}{4p_2d_1} + \dfrac{2d_1-p_1-2p_2}{4d_1}\right) &\geq 0 .\\
		\\
		\text{If we divide by }t_2^2 \text{ and distribute the 4 we have} \\
		\\
		\left(\frac{d_2^2}{d_1p_2} - \frac{t_1}{t_2}\right)^2 - \frac{d_2^2}{d_1p_2}\left(\frac{d_2^2}{p_2d_1} + \dfrac{2d_1 - p_1 - 2p_2}{d_1}\right) &\geq 0 .\\
	\end{align*}
	We now let $t = \frac{t_1}{t_2}$ and we determine when
	$$f(t) = \left(\frac{d_2^2}{d_1p_2} - t\right)^2 - \frac{d_2^2}{d_1p_2}\left(\frac{d_2^2}{p_2d_1} + \dfrac{2d_1 - p_1 - 2p_2}{d_1}\right) \geq 0$$
	since this will provide us with equivalent conditions to $c$ being real.
	\\
	\\Our $f(t)$ is a parabola in $t$ that is concave up with vertex $$\left(\frac{d_1^2}{d_2p_1} \text{ , } - \frac{d_2^2}{d_1p_2}\left(\frac{d_2^2}{p_2d_1} + \dfrac{2d_1 - p_1 - 2p_2}{d_1}\right)\right).$$
	Thus, the minimum of $f(t)$ is negative with one zero of $f(t)$ being greater than $\frac{d_1^2}{d_2p_1}$, so positive, and the other being less than $\frac{d_1^2}{d_2p_1}$, so unknown. We thus solve $f(t) = 0$ to determine the conditions for $f(t) \geq 0$, recalling that we really only care about $t < 0$ since $t_1 > 0$ and $t_2 < 0$.
	\begin{align*}
		f(t) &= 0 \\
		\left(\frac{d_2^2}{d_1p_2} - t\right)^2 - \frac{d_2^2}{d_1p_2}\left(\frac{d_2^2}{p_2d_1} + \dfrac{2d_1 - p_1 - 2p_2}{d_1}\right) &= 0 \\
		\left(\frac{d_2^2}{d_1p_2} - t\right)^2 &= \frac{d_2^2}{d_1p_2}\left(\frac{d_2^2}{p_2d_1} + \dfrac{2d_1 - p_1 - 2p_2}{d_1}\right) \\ 
		\frac{d_2^2}{d_1p_2} - t &= \pm \sqrt{\frac{d_2^2}{d_1p_2}\left(\frac{d_2^2}{p_2d_1} + \dfrac{2d_1 - p_1 - 2p_2}{d_1}\right)} \\
		t &= \frac{d_2^2}{d_1p_2} \pm \sqrt{\frac{d_2^2}{d_1p_2}\left(\frac{d_2^2}{p_2d_1} + \dfrac{2d_1 - p_1 - 2p_2}{d_1}\right)}  \\
	\end{align*}
	This provides us with $t$ such that $f(t) = 0$, but we drop the case with $+$ since we only want $t< 0$ and we observe when we have $t < 0$.
	\begin{align*}
		t &= \frac{d_2^2}{d_1p_2} - \sqrt{\frac{d_2^2}{d_1p_2}\left(\frac{d_2^2}{p_2d_1} + \dfrac{2d_1 - p_1 - 2p_2}{d_1}\right)} \\
		t&= \frac{d_2^2}{d_1p_2} - \frac{d_2}{d_1}\sqrt{\frac{d_2^2}{p_2^2} + \dfrac{2d_1 - p_1 - 2p_2}{p_2}} \leq 0
	\end{align*}
	where the final inequality is true by $\dfrac{2d_1 - p_1 - 2p_2}{p_2} \geq 0$. Observe that $t = 0$ only when $\dfrac{2d_1 - p_1 - 2p_2}{p_2} = 0$, and by Lemma \ref{pvalueinequalities} we know this only happens when $\p_1$ is a trivial representation. For that setting, please see Corollary \ref{ricequalscTtrivialp2}. In the current setting of $[\p_1, \p_1] \not\subset \h$, though, we may conclude that our $t$ for $f(t) = 0$ is negative.
	\\
	\\Since the $t$ found above where $f(t) = 0$ is negative and $f(t)$ is a parabola concave up, we can conclude that for $t \leq \frac{d_2^2}{d_1p_2} - \frac{d_2}{d_1}\sqrt{\frac{d_2^2}{p_2^2} + \dfrac{2d_1 - p_1 - 2p_2}{p_2}}$ we have $f(t) \geq 0$. Thus, having $(t_1, t_2)$ such that 
	\begin{align}
		\frac{t_1}{t_2} \leq \frac{d_2^2}{d_1p_2} - \frac{d_2}{d_1}\sqrt{\frac{d_2^2}{p_2^2} + \dfrac{2d_1 - p_1 - 2p_2}{p_2}} \label{tvals2}
	\end{align}
	is equivalent to having $c$ real and therefore $c > 0$ as well. This concludes our proof of \textbf{Claim 2}. 
	\\
	\\It will be helpful when we reach the end of the proof of \textbf{Claim 4} to go ahead and note that $\frac{t_1}{t_2} = t= \frac{d_2^2}{d_1p_2} - \frac{d_2}{d_1}\sqrt{\frac{d_2^2}{p_2^2} + \dfrac{2d_1 - p_1 - 2p_2}{p_2}}$ happens when $f(t) = 0$ which is when the discriminant of $c$ is 0 and $c_+ = c_-$.
	\\
	\\Now, we want to verify conditions for which not only $c>0$, but $\la = \frac{-2d_2}{p_2}\left(ct_2 + \frac{1}{2}\right) > 0$. 
	\\
	\\\textbf{\underline{Claim 3:}} For any $(t_1, t_2)$ satisfying (\ref{tvalsclaim}), there is one solution which can be obtained by $c_+$ and $\la = \frac{-2d_2}{p_2}\left(c_+t_2 + \frac{1}{2}\right)$ where $c_+$ takes the $+$ in (\ref{cvals})
	\\
	\\\underline{Proof of \textbf{Claim 3}:}
	Observe that $\la > 0$ occurs if and only if $ct_2 < \frac{-1}{2}$. Let us examine that inequality more closely:
	\begin{align*}
		ct_2 = \dfrac{-\left(\frac{d_2^2}{d_1p_2}t_2 - t_1\right) \pm \sqrt{\left(\frac{d_2^2}{d_1p_2}t_2 - t_1\right)^2 - 4\left(\frac{d_2^2}{d_1p_2}t_2^2\right)\left(\frac{d_2^2}{4p_2d_1} + \dfrac{2d_1-p_1-2p_2}{4d_1}\right)}}{2\frac{d_2^2}{d_1p_2}t_2} &< \frac{-1}{2}\\
		-\left(\frac{d_2^2}{d_1p_2}t_2 - t_1\right) \pm \sqrt{\left(\frac{d_2^2}{d_1p_2}t_2 - t_1\right)^2 - 4\left(\frac{d_2^2}{d_1p_2}t_2^2\right)\left(\frac{d_2^2}{4p_2d_1} + \dfrac{2d_1-p_1-2p_2}{4d_1}\right)} &> \frac{-d_2^2t_2}{d_1p_2} \\
		t_1 \pm \sqrt{\left(\frac{d_2^2}{d_1p_2}t_2 - t_1\right)^2 - 4\left(\frac{d_2^2}{d_1p_2}t_2^2\right)\left(\frac{d_2^2}{4p_2d_1} + \dfrac{2d_1-p_1-2p_2}{4d_1}\right)} &> 0 \tag{$\ast$}\label{ast1}
	\end{align*}
	This creates two cases for determining what conditions on $(t_1, t_2)$ provide $ct_2 < \frac{-1}{2}$, one in which we have the $+$ above in (\ref{ast1}) and another in which we have the $-$ above. Considering the case with the $+$, we recall that $t_1 > 0$, so as long as $c$ is real (and therefore positive), we have at least one solution. This completes the proof of \textbf{Claim 3}.
	\\
	\\We now turn our attention to the possibility of a second set of solutions determined by taking the $-$ in (\ref{ast1}). 
	\\
	\\\textbf{\underline{Claim 4:}} For any $(t_1, t_2)$ satisfying (\ref{tvalssecondsolution}), there is a second solution which can be obtained by $c_-$ and $\la = \frac{-2d_2}{p_2}\left(c_-t_2 + \frac{1}{2}\right)$ where $c_-$ takes the $-$ in (\ref{cvals})
	\\
	\\\underline{Proof of \textbf{Claim 4}:}
	In this case, we are looking for when the discriminant above in (\ref{ast1}) is less than $t_1^2$. Allow us to investigate:
	\begin{align*}
		\left(\frac{d_2^2}{d_1p_2}t_2 - t_1\right)^2 - 4\left(\frac{d_2^2}{d_1p_2}t_2^2\right)\left(\frac{d_2^2}{4p_2d_1} + \dfrac{2d_1-p_1-2p_2}{4d_1}\right) &< t_1^2 \\
		\frac{d_2^4}{d_1^2p_2^2}t_2^2 + t_1^2 - 2\frac{d_2^2}{d_1p_2}t_2t_1 - \frac{d_2^4}{d_1^2p_2^2}t_2^2 - \dfrac{d_2^2t_2^2(2d_1 - p_1 - 2p_2)}{d_1^2p_2} &< t_1^2 \\
		-2\frac{d_2^2}{d_1p_2}t_1t_2 - \dfrac{d_2^2t_2^2(2d_1 - p_1 - 2p_2)}{d_1^2p_2} &< 0 \\
		\text{now divide by } -t_2\frac{d_2^2}{d_1p_2}& \\
		2t_1 + t_2\dfrac{2d_1 - p_1 - 2p_2}{d_1} &< 0 \\
		\frac{t_1}{t_2} &> -\dfrac{2d_1 - p_1 - 2p_2}{2d_1}
	\end{align*}
	Also note that having a $(t_1, t_2)$ satisfy (\ref{ast1}) is equivalent to $ct_2 < \frac{-1}{2} < 0$, so we do not have to worry about if such a $(t_1, t_2)$ will not produce a $c > 0$ since we have $t_2 < 0$. Thus, for \begin{align}\label{secondsolution}
		-\dfrac{2d_1 - p_1 - 2p_2}{2d_1} < \frac{t_1}{t_2} < \frac{d_2^2}{d_1p_2} - \frac{d_2}{d_1}\sqrt{\frac{d_2^2}{p_2^2} + \dfrac{2d_1 - p_1 - 2p_2}{p_2}}
	\end{align}
	we have a second solution. Again, we remark that by Lemma \ref{pvalueinequalities}, when $\p_1$ is a trivial representation (the setting considered in Corollary \ref{ricequalscTtrivialp2}), we have the lower and upper bound being 0, and this is the only setting in which this can happen. We also exclude equality on the right as this happens if an only if the discriminant in c is 0 and in this case $c_+ = c_-$, so there is only one solution. This completes our proof of \textbf{Claim 4}.
	\\
	\\Furthermore, we observe that by Eqn. (\ref{c2.2}), we are able to determine what our metric is for a given solution, $T$ as $c$ is also known in terms of $T$.
	\\
	\\With the completion of \textbf{Claim 4}, we have our desired result providing necessary and sufficient conditions on $(t_1, t_2)$ for solutions to $ric = cT$. Moreover, we have provided the $c$ values dependent upon $T$, and the $(x_1, x_2)$ that determine our metric for the given $(t_1, t_2)$ and $c$ as well. \hfill $\blacksquare$
	\\
	\\
	\begin{corollary}
		Let $\g = \h \oplus \p_1 \oplus \p_2$ be noncompact simple with $[\p_1, \p_1] \subset \h$ and $\p_1$ not a trivial representation. For $G/H$ in this case, $ric = T$ has a solution $T$ if and only if $$t_1 = \frac{d_2^2}{d_1p_2}t_2^2 + \frac{d_2^2}{d_1p_2}t_2 + \dfrac{p_2(2d_1 - 2p_2) + d_2^2}{4d_1p_2} \text{ with } t_2 \in (-\infty, \frac{-1}{2}).$$
	\end{corollary}
	
	\vspace{.2in}
	
	\underline{Proof:}
	If $[\p_1, \p_1] \subset \h$, then we have $p_1 = 0$ in our $r_1$ and $r_2$ defined in Eqn.\ref{R1R2}. Thus, the solutions to $ric = T$ can be determined from Theorem \ref{ricequalsTgeneralp} to be any $(t_1, t_2)$ such that $t_1 =\frac{d_2^2}{d_1p_2}t_2^2 + \frac{d_2^2}{d_1p_2}t_2 + \dfrac{p_2(2d_1 - 2p_2) + d_2^2}{4d_1p_2}$ with  $t_2 \in (-\infty, \frac{-1}{2})$. \hfill $\blacksquare$
	\\
	\\
	
	\begin{corollary}
		Let $\g = \h \oplus \p_1 \oplus \p_2$ be noncompact simple with $[\p_1, \p_1] \subset \h$ and $\p_1$ not a trivial representation. For $G/H$ in this case, the equation $ric = cT$ has a solution if and only if $(t_1, t_2)$ is a pair satisfying $$\frac{t_1}{t_2} \leq \frac{d_2^2}{d_1p_2} - \frac{d_2}{d_1}\sqrt{\frac{d_2^2}{p_2^2} + \dfrac{2d_1 - 2p_2}{p_2}}$$
		where $t_1 < 0$ and $t_2 > 0$.
		\\
		\\When the above inequality is satisfied, there is always one solution which can be obtained by $c_+$ (where $c_+$ takes the $+$ in (\ref{cval2}) below) and $(x_1, x_2)$, the pair unique up to scaling given by $\frac{x_1}{x_2} = \frac{-2d_2}{p_2}(c_+t_2 + \frac{1}{2}).$
		\\In addition to this one solution, there is a second solution if and only if our pair $(t_1, t_2)$ with $t_1 > 0$ and $t_2 < 0$ satisfies 
		$$-\dfrac{d_1 - p_2}{d_1} < \frac{t_1}{t_2} < \frac{d_2^2}{d_1p_2} - \frac{d_2}{d_1}\sqrt{\frac{d_2^2}{p_2^2} + \dfrac{2d_1 - 2p_2}{p_2}}.$$
		The second solution is obtained by $c_-$ (where $c_-$ takes the $-$ in (\ref{cval2}) below) and $(x_1, x_2)$ the pair unique up to scaling given by $\frac{x_1}{x_2} = \frac{-2d_2}{p_2}(c_-t_2 + \frac{1}{2})$.
		\begin{align}\label{cval2}
			c = \dfrac{-\left(\frac{d_2^2}{d_1p_2}t_2 - t_1\right) \pm \sqrt{\left(\frac{d_2^2}{d_1p_2}t_2 - t_1\right)^2 - 4\left(\frac{d_2^2}{d_1p_2}t_2^2\right)\left(\frac{d_2^2}{4p_2d_1} + \dfrac{2d_1-2p_2}{4d_1}\right)}}{2\frac{d_2^2}{d_1p_2}t_2^2}
		\end{align}
		Moreover, $c_+ = c_-$ when the discriminant in $c$ is zero, and this happens precisely when $$\frac{t_1}{t_2} = \frac{d_2^2}{d_1p_2} - \frac{d_2}{d_1}\sqrt{\frac{d_2^2}{p_2^2} + \dfrac{2d_1 - 2p_2}{p_2}}.$$
	\end{corollary}
	\underline{Proof:}
	If $[\p_1, \p_1] \subset \h$ and $\p_1$ is not a trivial representation, then we have $p_1 = 0$ in our $r_1$ and $r_2$ defined in Eqn.\ref{R1R2}. Moreover, by Lemma \ref{pvalueinequalities} we have $\dfrac{2d_1-p_1-2p_2}{4d_1} = \dfrac{d_1-p_2}{2d_1} > 0$.
	Using Theorem \ref{ricequalscTgeneral} then we can say that $(t_1, t_2)$ satisfying
	\begin{align}\label{ast2}
		\frac{t_1}{t_2} \leq \frac{d_2^2}{d_1p_2} - \frac{d_2}{d_1}\sqrt{\frac{d_2^2}{p_2^2} + \dfrac{2d_1 - 2p_2}{p_2}}
	\end{align}
	is a sufficient and necessary condition to having a solution to $ric = cT$ obtained by $c_+$ from (\ref{cvaluesb}) below and $(x_1, x_2)$ satisfying $\frac{x_1}{x_2} = \frac{-2d_2}{p_2}\left(c_+t_2 + \frac{1}{2}\right)$.
	\\
	\\Moreover, having $(t_1, t_2)$ satisfy \[-\dfrac{d_1 - p_2}{d_1} < \frac{t_1}{t_2} < \frac{d_2^2}{d_1p_2} - \frac{d_2}{d_1}\sqrt{\frac{d_2^2}{p_2^2} + \dfrac{2d_1 - 2p_2}{p_2}} \tag{$\ast \ast$} \label{ast3}
	\] 
	is a sufficient and necessary condition to have a second solution, obtained by $c_-$ from (\ref{cvaluesb}) below and $(x_1, x_2)$ satisfying $\frac{x_1}{x_2} = \frac{-2d_2}{p_2}(c_- t_2 + \frac{1}{2})$.
	\begin{align}\label{cvaluesb}
		c = \dfrac{-\left(\frac{d_2^2}{d_1p_2}t_2 - t_1\right) \pm \sqrt{\left(\frac{d_2^2}{d_1p_2}t_2 - t_1\right)^2 - 4\left(\frac{d_2^2}{d_1p_2}t_2^2\right)\left(\frac{d_2^2}{4p_2d_1} + \dfrac{2d_1-2p_2}{4d_1}\right)}}{2\frac{d_2^2}{d_1p_2}t_2^2} 
	\end{align}
	Again, we do not have a second solution for $$\frac{t_1}{t_2} = \frac{d_2^2}{d_1p_2} - \frac{d_2}{d_1}\sqrt{\frac{d_2^2}{p_2^2} + \dfrac{2d_1 - 2p_2}{p_2}}$$
	as this is precisely when $c_+ = c_-$ since the discriminant of $c$ is 0.
	\hfill $\blacksquare$
	\\
	\\
	\begin{corollary}\label{ricequalsTtrivialad}
		Let $\g = \h \oplus \p_2 \oplus \p_1$ be noncompact simple with $\p_1$ trivial. For $G/H$ in this case, $ric = T$ has a solution $T$ if and only if $t_1 = d_2^2t_2^2 + d_2^2t_2 + \frac{d_2^2}{4}$ for $t_2 \in (-\infty, \frac{-1}{2})$.
	\end{corollary}
	\underline{Proof:}
	If $\p_1$ is a trivial representation then $p_1 = 0$ and $p_2 = d_1 = 1$ (See Lemma \ref{pvalueinequalities}) in our $r_1$ and $r_2$ defined in Eqn.\ref{R1R2}. Therefore, just as before, we use the result in Theorem \ref{ricequalsTgeneralp} to determine the solutions to $ric = T$. In this case, we have solutions of the form $t_1 = d_2^2t_2^2 + d_2^2t_2 + \frac{d_2^2}{4}$ with $t_2 \in (-\infty, \frac{-1}{2})$. \hfill $\blacksquare$
	\\
	\\
	\begin{corollary}\label{ricequalscTtrivialp2}
		Let $\g = \h \oplus \p_1 \oplus \p_2$ be noncompact simple with $\p_1$ trivial. For $G/H$ in this case, the equation $ric = cT$ has a solution for any given $(t_1, t_2)$ with $t_1 > 0$ and $t_2 < 0$. In this case, $c$ is defined by \ref{cval3} below and our inner product is defined by the $(x_1, x_2)$ pair unique up to scaling satisfying $\frac{x_1}{x_2} = -2d_1(ct_2 + \frac{1}{2})$.
		\begin{align}\label{cval3}
			c = \dfrac{-(d_2^2t_2 - t_1) + \sqrt{t_1^2 - 2d_2^2t_1t_2 }}{2d_2^2t_2^2}
		\end{align}
	\end{corollary}
	
	\vspace{.2in}
	
	\underline{Proof:}
	Once again we have $p_1 = 0$ and $p_2 = d_1 = 1$ (See Lemma \ref{pvalueinequalities}) in our $r_1$ and $r_2$ defined in Eqn. \ref{R1R2}. Moreover, by Lemma \ref{pvalueinequalities} we know that $\dfrac{2d_1-p_1-2p_2}{4d_1} =  0$, which we discussed beforehand in Remark \ref{specialpvalues} as being a situation which would present some changes in the types of solutions that arise in the proof of Theorem \ref{ricequalscTgeneral}. We discuss those changes now.
	\\
	\\Note that if we simply plug into the first solution to $ric = cT$ from Theorem \ref{ricequalscTgeneral}, then we have any $(t_1, t_2)$ such that $\frac{t_1}{t_2} \leq 0$ as solutions with the $c$ value being $c_+$ from
	\begin{align}\label{cvaluesa}
		c = \dfrac{-(d_2^2t_2 - t_1) \pm \sqrt{t_1^2 - 2d_2^2t_1t_2 }}{2d_2^2t_2^2}.
	\end{align} 
	\\
	\\Recalling the proof of \textbf{Claim 2} in Theorem \ref{ricequalscTgeneral}, we saw that $\dfrac{2d_1 - p_1 - 2p_2}{p_2} = 0$ causes us to have $c$ values that are real if and only if such that $\frac{t_1}{t_2} \leq 0$, and we saw in the proof of \textbf{Claim 3} that as long as $c$ was real, we had a solution with $c_+$ as the $c$ value. However, as was mentioned in the proof of \textbf{Claim 2}, we are interested only in $t_2 < 0$ and $t_1 > 0$ by necessity of $c , r_1 > 0$ and $r_2 < 0$. Thus, we must exclude the $\frac{t_1}{t_2} = 0$ from our solution set, implying that the sufficient and necessary condition to having a solution is $\frac{t_1}{t_2} < 0$. Therefore, for any $(t_1, t_2)$ with $t_1 > 0$ and $t_2<0$ gives a solution to $ric = cT$.
	\\
	\\Moreover, as the proof of \textbf{Claim 4} in Theorem \ref{ricequalscTgeneral} shows, $\dfrac{2d_1 - p_1 - 2p_2}{p_2} = 0$ in this case implies more than one solution does not exist because it would require $0< \frac{t_1}{t_2} < 0$. Thus, we do not have $c_-$ in this case. This gives us the desired result.  \hfill $\blacksquare$
	
	\newpage

	\section{5. $SO_0(1,7)/G_2$}\label{SO17}

	\vspace{.2in}

	\noindent In this section, at the Lie algebra level, we are working with $\so(1,7) = \g_2 \oplus \p_1 \oplus \p_2$ in which $\p_1 \simeq \p_2$ with dim$\p_1 = $ dim$\p_2 = 7$. The goal will be to determine the $ad_{\g_2}$ invariant $T(.,.)$ for which we have solutions to $ric = T$ and $ric = cT$, also supplying a way to find $c > 0$ in the second equation. 
	\\
	\\In the subsequent material, we will work with an orthonormal basis with respect to our fixed inner product $\langle . , . \rangle = B_{\p_2} - B_{\p_1}$ on $\p = \p_1 \oplus \p_2$, $\{x_1, ..., x_{14}\}$, with $\{x_1, ..., x_7\}$ a basis for $\p_1$ and $\{x_8, ..., x_{14}\}$ a basis for $\p_2$. As will be discussed in Lemma \ref{symmetricso17map} and Corollary \ref{TFormso17}, we are interested in working with a particular nice choice of an $\langle . ,. \rangle$ orthonormal basis that will make our $ad_{\g_2}$ equivariant maps have matrices where the blocks are diagonal. Corollary \ref{TFormso17}, in particular, show us that our $T$ is determined by $(t_1, t_2, t_3)$ such that $t_1 = T(x_1, x_1)$, $t_2 = T(x_8, x_8)$, and $t_3 = T(x_1, x_8)$ for our choice of basis. Therefore, we also have that $ric$ is determined by $(r_1, r_2, r_3)$ such that $r_1 = ric(x_1, x_1)$, $r_2 = ric(x_8, x_8)$, and $r_3 = ric(x_1, x_8)$. 
	\\
	\begin{remark}
		Our basis of choice is provided in Appendix \ref{AppA}. To find this basis, we used Sympy in Python (\cite{SymPy}) to determine our basis and compute $ric$ (\cite{code}). To get the diagonal blocks we wanted for our equivariant maps, all that was required was a re-ordering of the already obtained orthonormal basis with respect to $\langle . , . \rangle$.	
	\end{remark}
	
	\begin{remark}
		Due to our use of Mathematica (\cite{Mathematica}), there was no feasible way to determine the metric for a given solution $T$ to $ric = T$ or $ric = cT$. We are, however, able to determine our $c$ value in terms of $T$ in the latter equation. Previously, we were able to determine the metric as we could visibly see the substitution being made in the proof. Unfortunately, we were unable to see this kind of step with Mathematica.
	\end{remark}
	
	\vspace{.2in}
	
	\begin{theorem}\label{RiceqT}
		For $SO_0(1,7)/G_2$ with $(t_1, t_2, t_3)$ being defined as above, there is a $T$ such that $ric = T$ if and only if $(t_1, t_2, t_3)$ is such that:
		\\
		\\\textbf{Case 1} If $t_3 = 0$ then $(t_1, t_2, 0)$ is such that $t_1 = 6t_2^2 + 6t_2 + \frac{15}{8}$ with $t_2 < \frac{-1}{2}$.
		\\
		\\\textbf{Case 2} If $t_3 \not= 0$ then $(t_1, t_2, t_3)$ is such that $$t_1 = \begin{cases}
			f_1(t_2, t_3) \text{ , }& t_2 \leq \frac{-3}{4} \text{ and } |t_3| > 0 \\
			f_1(t_2, t_3) \text{ , }& -\frac{3}{4} < t_2 \leq \frac{-1}{2} \text{ and } |t_3| > 0 \text{ and } |t_3| \not=\frac{1}{2} \sqrt{\frac{3}{2}}\sqrt{4t_2 + 3} \\
			f_1(t_2, t_3) \text{ , }& t_2 > \frac{-1}{2} \text{ and } |t_3| >  \frac{\sqrt{3}}{4}\sqrt{2t_2  + 1}  \text{ and } |t_3| \not=\frac{1}{2} \sqrt{\frac{3}{2}}\sqrt{4t_2 + 3} \\
			\frac{3}{4} \text{ , }& t_2> \frac{-3}{4} \text{ and } |t_3| =\frac{1}{2} \sqrt{\frac{3}{2}}\sqrt{4t_2 + 3}
		\end{cases}$$
		where $f_1$ is described by Eqns.\ref{ricTsolutionso17} below.
	\end{theorem}
	
	\vspace{.2in}
	
	\begin{theorem}\label{RiceqcT}
		For $SO_0(1,7)/G_2$ with $(t_1, t_2, t_3)$ being defined as above, there is a $T$ such that $ric = cT$ for some $c>0$ if and only if $T$ is determined by $(t_1, t_2, t_3)$ such that:
		\\
		\\\textbf{Case 1} If $t_3 = 0$ then $(t_1, t_2, t_3)$ is such that $\frac{1}{3} \left(-\sqrt{5}-2\right) \leq \frac{t_2}{t_1} < 0$. In this case for a given solution $(t_1, t_2, t_3)$, $c = \frac{c_0}{t_1}$ where $c_0$ is the solution(s) to the implicit equation $c_0 = 6(c_0\frac{t_2}{t_1})^2 + 6c_0\frac{t_2}{t_1} + \frac{15}{8}$.
		\\
		\\\textbf{Case 2} If $t_3 \not= 0$ then $(t_1, t_2, t_3)$ is such that $\frac{t_2}{t_1} = l$ and $\frac{|t_3|}{t_1} = m$ in the 14 regions defined in Step 6 (See (\ref{Step 6})) below.
		Unless $(t_1, t_2, t_3)$ is a multiple of $(\frac{3}{4}, r_2, r_3)$ with $r_2 > \frac{-3}{4}$ and $|r_3| = \frac{1}{2}\sqrt{\frac{3}{2}}\sqrt{4r_2 + 3}$, for a given solution $(t_1, t_2, t_3)$, we have $c = \frac{c_0}{t_1}$ where $c_0$ is the solution(s) to the implicit equation $c_0 = f_1(c_0\dfrac{t_2}{t_1}, c_0\dfrac{t_3}{t_1})$ with $f_1$ is as defined in Theorem \ref{RiceqT}.  If $(t_1, t_2, t_3)$ is a multiple of $(\frac{3}{4}, r_2, r_3)$ with $r_2 > \frac{-3}{4}$ and $|r_3| = \frac{1}{2}\sqrt{\frac{3}{2}}\sqrt{4r_2 + 3}$, then $c = \frac{c_1}{t_1}$ where $c_1 = \frac{3}{4}$. 
	\end{theorem} 
	\vspace{.2in}
	\begin{remark}
		In what follows is a series of steps explaining the process. When we get to the end of Step 5 (see \ref{Step 5}) with the code reference provided, we have a proof for Theorem \ref{RiceqT} above. Similarly, when we get to the end of Step 6 (\ref{Step 6}), we have a proof of Theorem \ref{RiceqcT} above.
	\end{remark}
	
	\newpage
	
	\begin{landscape}
		\noindent Let $f_1(t_2, t_3)$ be the first root in the set of roots (in increasing order) to the following polynomial in $t_1$
		\\
		\begin{align}\label{ricTsolutionso17}
			128 t_1^3+t_1^2 \left(-768 t_2^2-768 t_2-432\right)+ t_1\left(1152 t_2^2+1536 t_2 t_3^2+1152 t_2+1536 t_3^2+432\right)-432 t_2^2-432 t_2-768 t_3^4-288 t_3^2-135
		\end{align}
		This root has the following radical form. This is not necessarily equavlent to $f_1(t_2, t_3)$ as the obtaining of radical forms for the roots of polynomials in which some coefficients are parameters can result in a loss of values.
		\begin{align*}
			&\frac{1}{8}\left(16 t_2^2+16 t_2+9\right)\\
			&+\frac{1}{8} \left(\sqrt[3]{\left(16 t_2^2+16 t_2+3\right)^3-192 (4 t_2+3) (2 t_2 (4 t_2+5)+5) t_3^2+8 \sqrt{2} \sqrt{-\left((4 t_2+1)^3-72 t_3^2\right) \left(3 (4 t_2+3)^2 t_3+16 t_3^3\right)^2}+1536 t_3^4} \right) \notag\\
			&+\frac{1}{8}\left(\frac{(4 t_2+1)^2 (4 t_2+3)^2}{\sqrt[3]{\left(16 t_2^2+16 t_2+3\right)^3-192 (4 t_2+3) (2 t_2 (4 t_2+5)+5) t_3^2+8 \sqrt{2} \sqrt{-\left((4 t_2+1)^3-72 t_3^2\right) \left(3 (4 t_2+3)^2 t_3+16 t_3^3\right)^2}+1536 t_3^4}}\right) \notag\\
			&-\frac{1}{8}\left(\frac{256 (t_2+1) t_3^2}{\sqrt[3]{\left(16 t_2^2+16 t_2+3\right)^3-192 (4 t_2+3) (2 t_2 (4 t_2+5)+5) t_3^2+8 \sqrt{2} \sqrt{-\left((4 t_2+1)^3-72 t_3^2\right) \left(3 (4 t_2+3)^2 t_3+16 t_3^3\right)^2}+1536 t_3^4}} \right) \notag
		\end{align*}
	\end{landscape}
	
	\includegraphics[width=\textwidth]{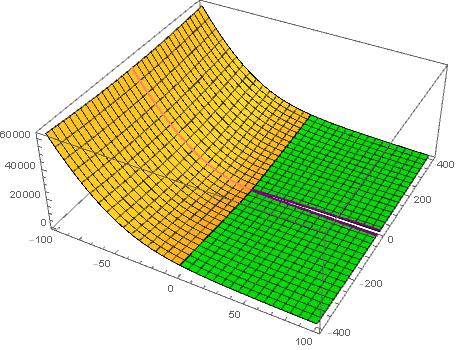}
	
	\noindent The above image was produced in Mathematica for the image of $ric$ in Theorem \ref{RiceqT}. In the below descriptions of parameters, $z$ is the third parameter describing $\Phi$, our positive definite equivariant map discussed in Step 1 (See \ref{Step 1}).
	\\
	\\The gold sheet curving up has a pink strip that can be faintly seen in the middle. The gold is the $z, t_3 \not= 0$ with $t_2 \leq \frac{-3}{4}$ solution but only graphed out to $-100 \leq t_2$, and the pink strip is $z, t_3 = 0$ which has $t_2 \leq \frac{-1}{2}$, also only graphed out to $-100 \leq t_2$. The gold and the pink were graphed with $0 < |t_3| \leq 400$. 
	\\
	\\The green is the image where $z \not= 0$, $t_2 > \frac{-1}{2}$, and  $|t_3| > \frac{1}{2}\sqrt{\frac{3}{2}}\sqrt{3 + 4t_2}$. Here, we only graphed out to $t_2 \leq 100$ with $|t_3| \leq 400$.
	\\
	\\There are other colors that can be seen in the image above describing the other parts of our solution to $ric = T$. To help show them somewhat more clearly, we have the following image which has  more restrictive bounds on the parameters ($t_2 \leq 50$ and $|t_3|< 50$ instead of the values used above of 100 and 400, respectively). Still, in this image, there are some strips that are hard to see (such as the solutions where $t_1 = \frac{3}{4}$). The piece going upward is the image where $-\frac{3}{4} < t_2 \leq -\frac{1}{2}$ and the the rest is a combination of the image where $-\frac{1}{2} < t_2$ and $\frac{\sqrt{3}}{4}\sqrt{2t_2  + 1} < |t_3| < \frac{1}{2} \sqrt{\frac{3}{2}}\sqrt{4t_2 + 3}$ and the image where $t_1 = \frac{3}{4}$.
	\\
	
	\includegraphics[width=\textwidth]{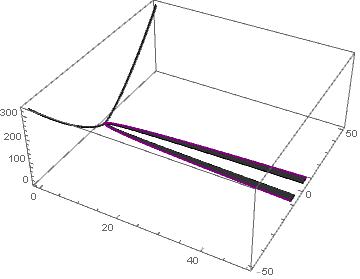}

	\newpage
	
	\noindent \underline{\textbf{An Overview of our Approach:}}
	\\
	\\In the case of $SO_0(1,7)/G_2$, we know from I.16 in \cite{Kerr} and our dual process in Theorem \ref{ClassificationTheorem} that $\so(1,7) = \g_2 \oplus \p_1 \oplus \p_2$ is such that $\p_1 \simeq \p_2$ of dimension 7. We are interested in solving $ric = T$ and $ric = cT$, so we begin with the primary aim of finding a forumula for $ric$ which is dependent upon any $ad_{\g_2}$ invariant inner product on $\p = \p_1 \oplus \p_2$. To do this, we use a description of all such inner products and our equation for $ric(.,.)$ in Eqn.\ref{ric}. Working out the solutions to $ric = T$ and $ric = cT$ turns out to be an involved problem even in this smaller dimensional example, so we provide a step-by-step overview of the process involved in reaching our solution before getting into the details.
	\begin{itemize}
		\item[\textbf{Step 1}] Determine a description of an arbitrary $ad_{\g_2}$ invariant inner product $(.,.) = \langle \Phi . , . \rangle = \langle \phi . , \phi . \rangle$ where $\Phi : \p \to \p$ is a positive definite $ad_{\g_2}$ equivariant map and $\phi : \p \to \p$ is such that $\phi^2 = \Phi$ and is also an $ad_{\g_2}$ equivariant map. In this step, we determine that $\phi$ is dependent upon three variables $(a,b,c)$. We also determine that $\Phi$ is determined by three variables which we label $(x,y,z)$ and we note the polynomial relationship between $(a,b,c)$ and $(x,y,z)$.
		\item[\textbf{Step 2}] Using the arbitrary inner product description of $(.,.) = \langle \phi . ,  \phi. \rangle$ (since defining an orthonormal basis on $\p$ for $(.,.)$ is easiest with this description), we find a formula for $ric(.,.)$ for an arbitrary $ad_{\g_2}$ invariant inner prodcut. The formula obtained depends only upon $\phi$, $\langle . , .\rangle$, and $\{x_i\}$ where $\{x_i\}$ is a $\langle . , . \rangle$-orthonormal basis on $\p$.
		\item[\textbf{Step 3}] Compute $ric(.,.)$ to obtain a function dependent only upon our $(a,b,c)$ defining $\phi$ using the description found in Step 2. For this, we use SymPy in Python (\cite{SymPy}) and acquire three terms $r_1$, $r_2$, and $r_3$ which are scale-invariant rational functions.
		\item[\textbf{Step 4}] Using polynomial relationships between the $(a,b,c)$ determining $\phi$ and the $(x,y,z)$ determining $\Phi$, we use an algebraic geometry tool known as an elimination ideal to get our $r_i$ in terms of $(x,y,z)$. This step was made possible by the built-in elimination ideal function in Mathmematica (\cite{Mathematica}), and cuts the degree of the polynomials (in both the numerator and denominator of our $r_i$) in half. 
		\item[\textbf{Step 5}] Using Mathematica (\cite{Mathematica}), we use built-in functions and utilize the scale-invariance of our $r_i$ to find all $(t_1, t_2, t_3)$ such that $(r_1, r_2, r_3) = (t_1, t_2, t_3)$ for some $(x, y, z)$ defining our $ad_{\g_2}$ invariant inner product, thus solving the problem of $ric = T$.  
		\item[\textbf{Step 6}] Using built-in functions from Mathematica (\cite{Mathematica}), and by projecting $(r_1, r_2, r_3)$ onto a plane, we find all $(t_1, t_2, t_3)$ such that there is a $c$ and an $(x,y,z)$ with $(r_1, r_2, r_3) = c(t_1, t_2, t_3)$. We also provide a description of the $c$ values needed for a given $(t_1, t_2, t_3)$.
	\end{itemize}
	\vspace{.2in}
	As we work through these six steps, we provide the mathematics here, and we explain some of the methods used in our code which can be found in \cite{code}.
	\\
	\begin{remark}
		Step 4 turned out to be a necessary step for us. The following two steps were attempted using the $r_i$ in terms of $(a,b,c)$, but the Mathematica functions utilized did not finish processing, even after more than 24 hours of run-time on the author's computer. After that attempt, there was an attempt to use the OU Supercomputing Center for Education \& Research (OSCER) at the University of Oklahoma (OU). Even this was unsuccessful in completing the computation of $ric = T$ and $ric = cT$ with the $r_i$ being provided in terms of $(a,b,c)$.
	\end{remark}
	\begin{remark}
		In Step 5, we take advantage of the scale-invariance. Using the scale invariance also turned out to be necessary as the functions in Mathematica spent hours running with no output.
	\end{remark}
	\begin{remark}
		In Mathematica, we used \textcolor{red}{\textbf{\textit{AbsoluteTiming}}} to find the run times for obtaining solutions to $ric = T$ and $ric = cT$ using the methodology discussed in Step 5 and Step 6, respectively. The run time for solving $ric = T$ was 126.031 seconds (so just over 2 minutes), and the run time for solving $ric = cT$ was 346.735 seconds (so just under 6 minutes). There was no use of OSCER in these run times, only the author's computer.
	\end{remark}
	
	\vspace{.3in}

	\subsection*{Step 1}\label{Step 1}

	Any $ad_{\g_2}$ invariant inner product $(.,.)$ on $\p$ can be determined by $(v,w) = \langle  \Phi v , w \rangle = \langle \phi v, \phi w \rangle$ where $\Phi = \phi^2$ is a positive definite $ad_{\g_2}$ equivariant map and $\phi$ is an invertible symmetric map. More than that, $\phi$ is also an $ad_{\g_2}$ equivariant map since $\Phi$ has a matrix representation that is diagonal for a basis of eigenvectors of $\p$ which implies that $\phi$ also has a matrix representation that is diagonal for the same basis of eigenvectors. Thus, $\Phi \circ ad_{\g_2} = ad_{\g_2} \circ \Phi$ will imply $\phi \circ ad_{\g_2} = ad_{\g_2} \circ \phi$. In this first step, we seek to understand the form of $\Phi$ and $\phi$ so that we can ultimately understand any $ad_{\g_2}$ invariant inner product in terms of $\Phi$ or $\phi$, and so that we can understand the algebraic relationship between $\Phi$ and $\phi$.
	\begin{lemma}\label{symmetricso17map}
		For $\so(1,7) = \g_2 \oplus \p_1 \oplus \p_2$, let $\p = \p_1 \oplus \p_2$. If $M : \p \to \p$ is a symmetric (with respect to $\langle . , . \rangle$) $ad_{\g_2}$ equivariant map then for a nice choice of basis, $M$ has a matrix of the form:
		$$M = \begin{bmatrix}
			\begin{BMAT}(b)[2pt, 5cm, 5cm]{c:c}{c:c}
				aId_{\p_1} & cId \\
				cId & bId_{\p_2}\\
			\end{BMAT}
		\end{bmatrix}$$  
		where $a,b,c \in\R$. Here, unless $c = 0$, the nice basis of choice is $\{x_1, ..., x_7, Lx_1, ... Lx_7\}$ where $\{x_i\}$ is an orthonormal basis with respect to $\langle ., . \rangle$ on $\p_1$ and $L: \p_1 \to \p_2$ is an $ad_{\g_2}$ intertwining map defined by $L = \proj_{\p_2} \circ M|_{\p_1}$. If $c = 0$ then the basis on $\p$ can be any basis orthonormal with respect to $\langle . , . \rangle$.
	\end{lemma}
	\underline{Proof:}
	First, if $M$ is a symmetric (with respect to $\langle . , . \rangle$) $ad_{\g_2}$ equivariant map then we know that in general, $$M = \begin{bmatrix}
		\begin{BMAT}(b)[2pt, 5cm, 5cm]{c:c}{c:c}
			A & L^t \\
			L & B\\
		\end{BMAT}
	\end{bmatrix}$$ where $A : \p_1 \to \p_1$ is symmetric, $B : \p_2 \to \p_2$ is symmetric, $L$ is defined as in the the statement of our Lemma, and $L^t$ is transposed with respect to $\langle . , . \rangle$ (along with the symmetry of $A$ and $B$). Since $\p_1$ and $\p_2$ are irreducible, and since $A$ and $B$ are symmetric, we know by Schur's Lemma that $A = a Id_{\p_1}$ and $B = b Id_{\p_2}$. Moreover, since dim$\p_1 =$ dim$\p_2 = 7$, we know that $\p_1$ and $\p_2$ are irreducible representations of real type (See \cite{BroTom}), meaning that $a,b \in \R$.
	\\
	\\Again by Schur's Lemma, $L :\p_1 \to \p_2$ an isomorphism or 0. If 0, then we are done and $c = 0$ in the statement of the Lemma. If $L$ is an isomorphism then we know, again by Schur's Lemma, that $LL^t = \la Id_{\p_1}$ with $\la \in \R$. Moreover, by choosing $\{x_1, ... , x_7, Lx_1, ..., Lx_7\}$ as a basis for $\p$ where $\{x_1, ..., x_7\}$ is an orthonormal basis for $\p_1$, we know that $L$ becomes a diagonal matrix. Thus, we know that $L = cId$ with $c \in \R$. \hfill $\blacksquare$
	\\
	\begin{remark}
		It is worth noting that in the case that $L: \p_1 \to \p_2$ is not the 0 map, we have that any other $ad_{\g_2}$ intertwining map $\p_1 \to \p_2$ is a multiple of $L$. Indeed, if $N: \p_1 \to \p_2$ was another $ad_{\g_2}$ intertwining map, then $N^{-1} L = \la Id_{\p_1}$ with $\la \in \R$. Thus, $N = \la L$. Originally it would appear that our basis was dependent upon a variety of choices for $L$, but this shows that there is only one choice for $L$ up to scaling.
		\\
	\end{remark}
	\begin{corollary}\label{TFormso17}
		For $\so(1,7) = \g_2 \oplus \p_1 \oplus \p_2$, let $\p = \p_1 \oplus \p_2$. If $T(.,.)$ is an $ad_{\g_2}$ invariant bilinear form then
		$$T(v,w) = \begin{cases}
			t_1\langle v , w \rangle \text{ , } v,w \in \p_1 \\
			t_2\langle v,w \rangle \text{ , } v,w \in \p_2 \\
			t_3\langle v,w \rangle \text{ , } v \in \p_1 \text{ , } w\in \p_2
		\end{cases}$$
		where $(t_1, t_2, t_3) \in \R^3$.
	\end{corollary}
	\underline{Proof:}
	If $T(.,.)$ is an $ad_{\g_2}$ bilinear form then, $T(v,w) = \langle M v, w \rangle$ where $v,w \in \p$ and $M$ is symmetric with respect to $\langle . , . \rangle$ and is an $ad_{\g_2}$ equivariant map. By Lemma \ref{symmetricso17map} we know that for a nice choice of basis, $M$ can be described by:
	$$M = \begin{bmatrix}
		\begin{BMAT}(b)[2pt, 5cm, 5cm]{c:c}{c:c}
			t_1Id_{\p_1} & t_3Id \\
			t_3Id & t_2Id_{\p_2}\\
		\end{BMAT}
	\end{bmatrix}$$
	for $t_1, t_2, t_3 \in \R$. Thus, we have that 	$$T(v,w) = \langle M v , w \rangle = \begin{cases}
		t_1\langle v , w \rangle \text{ , } v,w \in \p_1 \\
		t_2\langle v,w \rangle \text{ , } v,w \in \p_2 \\
		t_3\langle v,w \rangle \text{ , } v \in \p_1 \text{ , } w\in \p_2,
	\end{cases}$$
	as desired. \hfill $\blacksquare$
	\\
	\begin{lemma} \label{so17intertwiningmaps}
		For $\so(1,7) = \g_2 \oplus \p_1 \oplus \p_2$, let $\p = \p_1 \oplus \p_2$. Any $ad_{\g_2}$ invariant inner product $(.,.)$ can be written as $(.,.) = \langle \Phi . , . \rangle = \langle \phi . , \phi . \rangle$ where $\Phi$ is positive definite, $\phi$ is symmetric and invertible, and both $\Phi$ and $\phi$ are $ad_{\g_2}$ equivariant maps. Moreover, for a nice choice of basis (See Lemma \ref{symmetricso17map}) we have the following descriptions of $\Phi$ and $\phi$: 
		
		$$\Phi = \begin{bmatrix}
			\begin{BMAT}(b)[2pt, 5cm, 5cm]{c:c}{c:c}
				xId_{\p_1} & zId \\
				zId & yId_{\p_2}\\
			\end{BMAT}
		\end{bmatrix} \text{ with } x,y > 0 \text{ and } xy - z^2 > 0$$ 
		
		$$\phi = \begin{bmatrix}
			\begin{BMAT}(b)[2pt, 5cm, 5cm]{c:c}{c:c}
				aId_{\p_1} & cId \\
				cId & bId_{\p_2}\\
			\end{BMAT}
		\end{bmatrix} \text{ with } ab - c^2 \not= 0.$$
	\end{lemma}
	\underline{Proof:}
	In the following proof, we first prove the statement regarding $\Phi$. Then, using the positive definiteness of $\Phi$ and properties of square roots, we then prove the statement regarding $\phi$. When proving the statement regarding $\phi$, we prove the statement in the $2\times 2$ matrice setting and then show that the statement holds in the $14 \times 14$ setting we are in. 
	\\
	\\We know that any $ad_{\g_2}$ invariant inner product on $\p = \p_1 \oplus \p_2$ is of the form $(.,.) = \langle \Phi . , . \rangle$ where $\Phi: \p \to \p$ is a symmetric (with respect to $\langle ., . \rangle$) positive definite $ad_{\g_2}$ equivariant map. Moreover, since $\Phi$ is symmetric with respect to $\langle .,. \rangle$, we know by Lemma \ref{symmetricso17map} that for a nice choice of basis we have $\Phi$ of the form 
	\setlength{\abovedisplayskip}{0pt} 
	\setlength{\belowdisplayskip}{0pt} 
	\setlength{\abovedisplayshortskip}{0pt} 
	\setlength{\belowdisplayshortskip}{0pt}
	$$\Phi = \begin{bmatrix}
		\begin{BMAT}(b)[2pt, 5cm, 5cm]{c:c}{c:c}
			xId_{\p_1} & zId \\
			zId & yId_{\p_2}\\
		\end{BMAT}
	\end{bmatrix}.$$
	Since $\Phi$ is positive definite, $tr\Phi = 7(x+y) > 0$, and $det\Phi = (xy - z^2)^7 > 0$, we can conclude that $x,y > 0$ and $xy - z^2 > 0$. Thus, $$P = \left\{\begin{bmatrix}
		\begin{BMAT}(b)[2pt, 5cm, 5cm]{c:c}{c:c}
			xId_{\p_1} & zId \\
			zId & yId_{\p_2}\\
		\end{BMAT}
	\end{bmatrix}: x,y > 0 \text{ and } xy - z^2 > 0\right\}$$ is the set describing all positive definite matrices on $\p$. This concludes the proof for the first statement regarding $\Phi$.
	\\
	\\By being positive definite, $\Phi$ has a square root matrix $\phi$ such that $\phi^2 = \Phi$. To complete the proof, we consider the $2\times 2$ matrix setting, determining what the collection of such $\phi$ is, and then we show that understanding the $2\times 2$ matrix setting is enough to determine the $14\times 14$ matrix setting our problem is placed within.
	\\
	\\In the $2\times 2$ setting, we consider positive definite $\Phi = \begin{bmatrix}
		x & z \\
		z & y
	\end{bmatrix}$. By being positive definite, there is a $C$ such that $L = \begin{bmatrix}
		\la_1 & 0 \\
		0 & \la_2 
	\end{bmatrix} = C \Phi C^{-1}$ where $\la_1, \la_2$ are the eigenvalues of $\Phi$. Since $trL = tr\Phi$ and $detL = det\Phi$, we know that $\Phi$ is positive definite if and only if $tr \Phi > 0$ and $det \Phi > 0$. Therefore, $\Phi$ is positive definite if and only if $x,y > 0$ and $xy - z^2 > 0$.
	\\Now $L$ has square root matrices described by, $l = \begin{bmatrix}
		\pm \sqrt{\la_1} & 0 \\
		0 & \pm \sqrt{\la_2}
	\end{bmatrix}$, and one can observe that all symmetric $\phi$ such that $\phi^2 = \Phi$ are such that $l = C \phi C^{-1}$ for some $l$. Since $l$ can be any diagonal matrix with nonzero determinant (because the $\la_i$ can be anything positive), it is then the case that $\phi$ can be any symmetric matrix with nonzero determinant. 
	\\
	\\Now we show that the collection of invertible, symmetric $2\times 2$ matrices generates all positive definite $2\times 2$ matrices by squaring, and then show how the $14 \times 14$ matrice setting has the same multiplication structure, completing our proof.
	\\
	\\Let $\phi = \begin{bmatrix}
		a & c \\
		c & b
	\end{bmatrix}$ with $ab - c^2 \not = 0$. In this case, $\phi^2 = \begin{bmatrix}
		a^2 + c^2 & c(a+b) \\
		c(a+b) & b^2 + c^2
	\end{bmatrix}$, which is positive definite if and only if $tr\phi^2 > 0$ and $det\phi^2 > 0$, and this is the case if and only if
	\setlength{\abovedisplayskip}{30pt} 
	\setlength{\belowdisplayskip}{30pt} 
	\setlength{\abovedisplayshortskip}{30pt} 
	\setlength{\belowdisplayshortskip}{30pt}
	\begin{align*}
		(a^2 + c^2)(b^2 + c^2) - c^2(a+b)^2 &> 0. \\
		\\
		\text{Now, }(a^2 + c^2)(b^2 + c^2) - c^2(a+b)^2 &= a^2b^2 + c^4 - c^2(2ab) \\
		&= (ab - c^2)^2, 
	\end{align*}
	and $(ab - c^2)^2 > 0$ if and only if $ab - c^2 \not = 0$. Thus, $\phi^2$ is clearly positive definite. Lastly, observe that $a^2 + c^2$ and $b^2 + c^2$ can be any positive number while $c(a+b)$ can be any non-negative number. Thus, in the $2\times 2$ setting, the invertible, symmetric $\{\phi\}$ generates the positive definite $\{\Phi\}$ by squaring each $\phi$.
	\\
	\\To extend this to the $14 \times 14$ matrix case is now a trivial observation that if, in the $2 \times 2$ setting, $\{\phi\}$ generates all possible positive definite matrices $\{\Phi\}$ by squaring each $\phi$, then by having the same multiplication structure: 
	
	\begin{align*}
		\begin{bmatrix}
			\begin{BMAT}(b)[2pt, 5cm, 5cm]{c:c}{c:c}
				aId_{\p_1} & cId \\
				cId & bId_{\p_2}\\
			\end{BMAT}
		\end{bmatrix}\begin{bmatrix}
			\begin{BMAT}(b)[2pt, 5cm, 5cm]{c:c}{c:c}
				aId_{\p_1} & cId \\
				cId & bId_{\p_2}\\
			\end{BMAT}
		\end{bmatrix} \\
		=\begin{bmatrix}
			\begin{BMAT}(b)[2pt, 5cm, 5cm]{c:c}{c:c}
				(a^2 + c^2)Id_{\p_1} & c(a+b)Id \\
				c(a+b)Id & (b^2 + c^2)Id_{\p_2}\\
			\end{BMAT}
		\end{bmatrix}
	\end{align*}
	we can see that the analogous statement is true in the $14\times 14$ case. That is, $P$ as defined above, is generated by 
	\setlength{\abovedisplayskip}{0pt} 
	\setlength{\belowdisplayskip}{0pt} 
	\setlength{\abovedisplayshortskip}{0pt} 
	\setlength{\belowdisplayshortskip}{0pt}
	$$\left\{\begin{bmatrix}
		\begin{BMAT}(b)[2pt, 5cm, 5cm]{c:c}{c:c}
			aId_{\p_1} & cId \\
			cId & bId_{\p_2}\\
		\end{BMAT}
	\end{bmatrix}: ab - c^2\not= 0 \right\}$$
	by squaring the matrices. Thus, any inner product $(.,.) = \langle \Phi . , .\rangle = \langle \phi^2 . , . \rangle = \langle \phi . , \phi . \rangle$, with $\Phi$ and $\phi$ having the form desired. Moreover, since $\phi$ has the same block form as $\Phi$, it is clear that $\Phi$ being $ad_{\g_2}$ equivariant implies that $\phi$ is $ad_{\g_2}$ equivariant. \hfill $\blacksquare$
	\setlength{\abovedisplayskip}{0pt} 
	\setlength{\belowdisplayskip}{0pt} 
	\setlength{\abovedisplayshortskip}{0pt} 
	\setlength{\belowdisplayshortskip}{0pt}
	\begin{remark}
		As was noted in our list of steps above, it will later be useful to note the polynomial relationship between the $(a,b,c)$ defining $\phi$ and the $(x,y,z)$ defining $\Phi$. The following can be observed from how our $\{\phi\}$ generates our $\{\Phi\}$ in the preceding Lemma:
		\setlength{\abovedisplayskip}{0pt} 
		\setlength{\belowdisplayskip}{0pt} 
		\setlength{\abovedisplayshortskip}{0pt} 
		\setlength{\belowdisplayshortskip}{0pt}
		\begin{align}\label{polynomial relationships}
			a^2 + c^2 &= x \notag\\
			b^2 + c^2 &= y \\
			c(a+b) &= z. \notag
		\end{align}
	\end{remark}
	\vspace{.3in}
	\subsection*{Step 2}\label{Step 2}
	In the present step, we will be finding a formula for $ric(.,.)$ that is dependent only upon $\phi$, our fixed inner product $\langle ., . \rangle$, and a $\langle . , . \rangle$ orthonormal basis, $\{x_i\}$. In the end, we provide a formula for $ric(x,y)$ in which one can determine how to work with matrix representations of $\phi$ and $ad_\p$ with respect to the basis $\{x_i\}$ (except for the Killing form term which we intentionally leave as is).
	\\
	\begin{lemma}
		Consider $G/H$ to be any homogeneous space for which $G$ is unimodular with reductive decomposition $\g = \h \oplus \p$. For an arbitrary $ad_\h$ invariant inner product $(.,.) = \langle \phi . , \phi .\rangle$, we have the following formula for $ric(.,.)$ in terms of the base inner product $\langle . , . \rangle$, a $\langle . , . \rangle$ orthonormal basis $\{x_i\}$, and $\phi$.
		$$ric(x,y)= \frac{-1}{2}\sum_{x_i} \langle \phi (\sum_{k}(\phi^{-1})_{ki} ad_\p((x_k))(x)) ,\phi (\sum_{k}(\phi^{-1})_{ki} ad_\p((x_k))(y))\rangle - \frac{1}{2}B(x,y)$$
		$$+ \frac{1}{4}\sum_{x_i,x_j} \langle\phi (\sum_{k}(\phi^{-1})_{ki} ad_\p((x_k))(\phi^{-1}(x_j)))  , \phi(x)\rangle\langle \phi(\sum_{k}(\phi^{-1})_{ki} ad_\p((x_k))(\phi^{-1}(x_j))) , \phi(y)\rangle$$ 
	\end{lemma}
	\underline{Proof:} Using Eqn.\ref{ric}, we have the following in which, by an abuse of notation, we consider $ad_\p(x)$ to be defined by $[x, .]_\p$ where $x \in \p$:
	
	\begin{align*}
		ric(x,y) &= \frac{-1}{2}\sum_{x_i} ([\phi^{-1}(x_i) , x]_\p , [\phi^{-1}(x_i), y]_\p) - \frac{1}{2}B(x,y) \\
		&+ \frac{1}{4}\sum_{x_i,x_j} ([\phi^{-1}(x_i), \phi^{-1}(x_j)] , x)([\phi^{-1}(x_i), \phi^{-1}(x_j)] , y) \\
		&= \frac{-1}{2}\sum_{x_i} (ad_\p(\phi^{-1}(x_i)(x)) , ad_\p(\phi^{-1}(x_i)(y))) - \frac{1}{2}B(x,y) \\
		&+ \frac{1}{4}\sum_{x_i,x_j} (ad_\p(\phi^{-1}(x_i))(\phi^{-1}(x_j)) , x)(ad_\p(\phi^{-1}(x_i))(\phi^{-1}(x_j)) , y) \\
		& = \frac{-1}{2}\sum_{x_i} \langle \phi \circ ad_\p(\phi^{-1}(x_i)(x)) ,\phi \circ ad_\p(\phi^{-1}(x_i)(y))\rangle - \frac{1}{2}B(x,y) \\
		&+ \frac{1}{4}\sum_{x_i,x_j} \langle\phi \circ ad_\p(\phi^{-1}(x_i))(\phi^{-1}(x_j)) , \phi( x)\rangle\langle \phi \circ ad_\p(\phi^{-1}(x_i))(\phi^{-1}(x_j)) , \phi(y)\rangle.
	\end{align*}
	\vfill
	\newpage
	\noindent Using the linearity of  $ad$ and the fact that $\phi^{-1}(x_i) = (\phi^{-1})_{1i}x_1 + ... + (\phi^{-1})_{ni}x_n$ where $(\phi^{-1})_{ij} = \langle \phi^{-1}x_j , x_i \rangle$ (i.e the matrix entries given by a matrix representation of $\phi^{-1}$ with respect to the basis $\{x_i\}$), we have our result:
	\begin{align}\label{ricphi}
		&ric(x,y) = \frac{-1}{2}\sum_{x_i} \langle \phi (\sum_{k}(\phi^{-1})_{ki} ad_\p((x_k))(x)) ,\phi (\sum_{k}(\phi^{-1})_{ki} ad_\p((x_k))(y))\rangle \\
		&- \frac{1}{2}B(x,y) + \notag\\
		&\frac{1}{4}\sum_{x_i,x_j} \langle\phi (\sum_{k}(\phi^{-1})_{ki} ad_\p((x_k))(\phi^{-1}(x_j)))  , \phi(x)\rangle\langle \phi(\sum_{k}(\phi^{-1})_{ki} ad_\p((x_k))(\phi^{-1}(x_j))) , \phi(y)\rangle \notag
	\end{align} \hfill $\blacksquare$
	\vspace{.3in}
	\subsection*{Step 3}\label{Step 3}
	Having our description of $ric(.,.)$ for an arbitrary inner product, we can see that it is described totally in terms of the Killing form $B(.,.)$, the background inner product $\langle . , . \rangle$, a basis $\{x_i\}$ on $\p$ for which $\langle . , . \rangle$ is orthonormal, the $ad_\p(x_k)$ maps, and an equivariant map $\phi$. We now turn our attention to finding $ric(x,y)$ for an arbitrary inner product on $\p = \p_1 \oplus\p_2$ in $\so(1,7) = \g_2 \oplus \p_1 \oplus \p_2$. Using Eqn.\ref{ricphi}, our process is as follows:
	\begin{enumerate}
		\item[i.] Fix a background inner product and get an orthonormal basis for $\p = \p_1 \oplus \p_2$.
		\item[ii.] Find the equivariant maps $\phi$ described in Lemma \ref{so17intertwiningmaps} by finding an appropriately ordered basis.
		\item[iii.] Find the matrices $ad_\p(x_k)$ for each $x_k$ coming from the basis in i.
		\item[iv.] Determine $ric(.,.)$ in terms of the $a,b,c$ defining $\phi$.
	\end{enumerate}

	\noindent To do this, we use Python (\cite{SymPy}) and Maple (\cite{Maple}) with the code found in \cite{code}. 
	\\
	\\To see the full details of how we developed our expressions for $ric$ can be found in \cite{mythesis} or can be followed in the code provided in \cite{code}.
	\\
	\\Once we run our code to produce $ric(x,y)$ (by using $\phi^{-1}$ instead of $\phi$ as our equivariant map to make terms simpler), we ran checks to ensure that we have what we expected. From Corollary \ref{TFormso17}, we expect for $ric(x_i, x_i) \not= 0$ for $x \leq 14$ and for $ric(x_{x_i}, x_{i+7}) = ric(x_{x_{i+7}}, x_{x_i}) \not = 0$ for $i \leq 7$. We further expect for all other terms to be zero. From more general properties of $ric$, we expect for $ric(x_i, x_j) = ric(x_j, x_i)$, and for the expressions defining $ric(x_i, x_j)$ to be scale invariant with respect to a multiple of $\phi$ as this describes our metric. Lastly, again from Corollary \ref{TFormso17}, we expect for $ric(x_1, x_1) = ... = ric(x_7, x_7)$, $ric(x_8, x_8) = ... = ric(x_{14}, x_{14})$, and $ric(x_1, x_8) = ... = ric(x_8, x_{14})$. In the end, we label $r_1 = ric(x_1, x_1)$, $r_2 = ric(x_8, x_8)$ and $r_3 = ric(x_1, x_8)$. The code for these checks is provided in \cite{code}. We provide the $r_1, r_2$, and $r_3$ in terms of $(a,b,c)$ below.
	\vfill
	\pagebreak
	\begin{landscape}
		\noindent$r_1 = \frac{9 a^4 \left(b^4+c^4\right)-36 a^3 b c^2 \left(b^2-c^2\right)+6 a^2 c^2 \left(b^4+20 b^2 c^2+c^4\right)+12 a b c^2 \left(b^4+5 b^2 c^2-2 c^4\right)+b^8+10 b^6 c^2+27 b^4 c^4+10 b^2 c^6+10 c^8}{24 \left(c^2-a b\right)^4}$
		\\
		\\
		$r_2 = \frac{-2 \left(\left(a^2+c^2\right)^2+2 c^2 (a+b)^2+\left(b^2+c^2\right)^2\right) \left(c^2-a b\right)^2+c^2 \left(2 a^2 b+a \left(b^2+c^2\right)+b^3+3 b c^2\right)^2+c^2 (3 a+b)^2 \left(a^2+c^2\right)^2+2 \left(a^3 b+2 a^2 c^2+3 a b c^2+b^2 c^2+c^4\right)^2-12 \left(c^2-a b\right)^4}{24 \left(c^2-a b\right)^4}$
		\\
		\\
		$r_3 = -\frac{c (a+b) \left(b^2+c^2\right) \left(3 a^2+b^2+4 c^2\right)^2}{24 \left(c^2-a b\right)^4}$
	\end{landscape}

	\vspace{.3in}
	
	\subsection*{Step 4}\label{Step 4}

	As we observed after the proof of Lemma \ref{so17intertwiningmaps} with Eqns.\ref{polynomial relationships}, the relationship between $\phi$ and $\Phi$ provides us with polynomial relationships between the $(a,b,c)$ defining $\phi$ and the $(x,y,z)$ defining $\Phi$:
	\begin{align*}
		a^2 + c^2 &= x \\
		b^2 + c^2 &= y \notag\\
		c(a+b) &= z. \notag
	\end{align*}
	Moreover, we have that $det\Phi = (det\phi)^2$. Using the polynomial relationships and recognizing that the denominator of each $r_i$ is a constant multiple of $(det\phi)^4$, we are able to use a function in Mathematica, \textcolor{red}{\textbf{\textit{Eliminate}}}, that finds elimination ideals (see the first three chapters of \cite{Cox} for more information on elimination ideals) to determine what the numerator of each $r_i$ is in terms of $(x,y,z)$. Checking that our newfound $r_i$ terms are correct is done by substituting back in for $x$, $y$, and $z$. The Mathematica code for the $r_1, r_2$, and $r_3$ in terms of $(a,b,c)$, the process of finding the $r_1, r_2$, and $r_3$ in terms of $(x,y,z)$, and the checks to ensure that our new $r_i$ terms are correct can be found in \cite{code}.
	\\
	\\Below, we provide our rational functions describing $r_1, r_2$, and $r_3$ in terms of $(x,y,z)$.
	\setlength{\abovedisplayskip}{0pt} 
	\setlength{\belowdisplayskip}{0pt} 
	\setlength{\abovedisplayshortskip}{0pt} 
	\setlength{\belowdisplayshortskip}{0pt}
	\begin{align*}
		r_1 &= \frac{9 x^2 y^2-18 x y z^2+y^4+6 y^2 z^2+18 z^4}{24 \left(z^2-x y\right)^2} \\
		r_2 &= \frac{-3 x^2 \left(4 y^2-3 z^2\right)-2 x \left(y^3-12 y z^2\right)+3 y^2 z^2-6 z^4}{24 \left(z^2-x y\right)^2} \\
		r_3 &= -\frac{y z (3 x+y)^2}{24 \left(z^2-x y\right)^2}
	\end{align*}
	\vspace{.3in}
	
	\subsection*{Step 5}\label{Step 5}
	The goal now is to find the solutions to $ric = T$ which we do by finding equations describing the $(t_1, t_2, t_3)$ such that there exists an $(x,y,z)$ with $x,y > 0$ and $xy - z^2 > 0$ with $(r_1, r_2, r_3) = (t_1, t_2, t_3)$. As can be seen in the $r_i$ provided in Step 4, if $(x,y,z) \mapsto (r_1, r_2, r_3)$, then $(x,y,-z) \mapsto (r_1, r_2, -r_3)$. Moreover, $r_3 = 0$ if and only if $z = 0$. This means that if we find solutions in the case of $z > 0$ then we have solutions to case when $z < 0$ given by $(t_1, t_2, -t_3)$ where $(t_1, t_2, t_3)$ is a solution to the $z > 0$ case. Therefore, we may consider solutions for when $t_3 > 0$ and obtain solutions for when $t_3 \not = 0$ by using $|t_3|$ in our conditions on $t_3$.
	\\
	\\
	Moreover, if $z = 0$ then we already have our solution since this is the case provided by Theorem \ref{ricequalsTgeneralp}. However, to show how our methods of using Mathematica cohere with that solution, we also provide the solution produced by Mathematica in the z = 0 setting in \cite{code}. We also provide in \cite{code} the solution to $ric = cT$ in the setting with $z = 0$ to show how the Mathematica code coheres with Theorem \ref{ricequalscTgeneral} as well. 
	\\
	\\Now we focus our attention on $(r_1, r_2, r_3) = (t_1, t_2, t_3)$, but with $t_3, r_3, z > 0$. Finding the solutions to $ric = T$ in this case requires us to use a combination of functions in Mathematica: \textcolor{red}{\textbf{\textit{Resolve}}} and \textcolor{red}{\textbf{\textit{Exists}}}. We sought to utilize these functions with $r_1, r_2$, and $r_3$ as is without using scale invariance; however, after hours of running, there was no output. Thus, we utilized the scale invariance and set $det\Phi = xy - z^2 = 1$. With $det\Phi = 1$, $r_1$ and $r_2$ became polynomials in terms of $(x,y)$ with rational coefficients, but we could not eliminate the $z$ in $r_3$, so $r_3$ became a polynomial in terms of $(x,y,z)$ with rational coefficients. Since $z$ was not eliminated completely, when using the combination of \textcolor{red}{\textbf{\textit{Resolve}}} and \textcolor{red}{\textbf{\textit{Exists}}} in Mathematica, we had to ensure that the polynomial relationship $xy - z^2 = 1$ was provided as a constraint. The code for this is provided in \cite{code}. We provide the new $r_1, r_2$, and $r_3$ below:
	\setlength{\abovedisplayskip}{10pt} 
	\setlength{\belowdisplayskip}{10pt} 
	\setlength{\abovedisplayshortskip}{10pt} 
	\setlength{\belowdisplayshortskip}{10pt}
	\begin{align*}
		r_1 &= \frac{1}{24} \left(9 x^2 y^2+6 x \left(y^2-3\right) y+y^4-6 y^2+18\right) \\
		r_2 &= \frac{1}{24} (9 x^3 y + x y (-12 + y^2) - 3 (2 + y^2) + x^2 (-9 + 6 y^2)) \\
		r_3 &= -\frac{1}{24}y (3 x + y)^2 z \text{ and } xy - z^2 = 1
	\end{align*}
	
	\begin{remark}
		One might think that eliminating $z$ could be done using elimination theory from Algebraic Geometry (See Chapters 1-3 of \cite{Cox} for this) as we did to get $r_1, r_2$, and $r_3$ in terms of $(x,y,z)$. However, the \textcolor{red}{\textbf{\textit{Eliminate}}} function, through its use of Grobner bases, only guarantees (without additional assumptions and algorithms being taken into consideration) an algebraic closure of the image of the function in consideration. Thus, when seeking to eliminate $z$ in $r_3$ by using $xy - z^2 = 1$, the resulting equation for $r_3$ in terms of $(x,y)$ need not have the same graph as $r_3$ in terms of $(x,y,z)$. Instead, we expect the algebraic closure of the graph of $r_3$. To find the image itself is called (in \cite{Cox}) the \textbf{implicitization problem} for functions that are rational or polynomial, and is much more involved than just elimination.
		\\
	\end{remark}
	\noindent The output that Mathematica provides from \textcolor{red}{\textbf{\textit{Resolve}}} and \textcolor{red}{\textbf{\textit{Exists}}} describes the region that is the image of $ric$ strictly in terms of $(t_1, t_2, t_3)$. The output was given by multiple regions with roots to polynomials being used to describe the region. However, we were able to simplify the regions and we were able to use the \textcolor{red}{\textbf{\textit{ToRadicals}}} function to get the roots of polynomials as functions written explicitly in terms of two variables. To simplify the regions, one observation we made was that the roots of the polynomials given turned out to be described by the same function in each region. Once we realized that, we were able to simplify the presentation of the regions, and the simplified form is given below, describing our $(t_1, t_2, t_3)$ such that $(r_1, r_2, r_3) = (t_1, t_2, t_3)$ for some $(x,y,z)$ with $z > 0$:
	$$t_1 = \begin{cases}
		f_1(t_2, t_3) \text{ , }& t_2 \leq \frac{-3}{4} \text{ and } t_3 > 0 \\
		f_1(t_2, t_3) \text{ , }& -\frac{3}{4} < t_2 \leq \frac{-1}{2} \text{ and } t_3 > 0 \text{ and } t_3 \not=\frac{1}{2} \sqrt{\frac{3}{2}}\sqrt{4t_2 + 3} \\
		f_1(t_2, t_3) \text{ , }& t_2 > \frac{-1}{2} \text{ and } t_3 >  \frac{\sqrt{3}}{4}\sqrt{2t_2  + 1}  \text{ and } t_3 \not=\frac{1}{2} \sqrt{\frac{3}{2}}\sqrt{4t_2 + 3} \\
		\frac{3}{4} \text{ , }& t_2> \frac{-3}{4} \text{ and } t_3 =\frac{1}{2} \sqrt{\frac{3}{2}}\sqrt{4t_2 + 3}
	\end{cases}$$
	where $f_1(t_2, t_3)$ is described in Eqn.\ref{ricTsolutionso17}. 
	\\
	\\This provides us with the solution desired for $z, t_3 > 0$ which finishes \textbf{Case 2} of Theorem \ref{RiceqT} regarding $t_3 \not= 0$.
	\\
	\\As mentioned before, we provide the Mathematica code for \textbf{Case 1} of Theorem \ref{RiceqT} in \cite{code}. The solution in this case follows the same program as \textbf{Case 2}, but with $z = r_3 = 0$. With those conditions, we use \textcolor{red}{\textbf{\textit{Resolve}}} and \textcolor{red}{\textbf{\textit{Exists}}} to find the following solution:
	$$t_1 = \frac{3}{8}(5 + 16t_2 + 16t_2^2) \text{ for } t_2 < -\frac{1}{2}$$
	
	\underline{Proof of Theorem \ref{RiceqT}}: Steps 1 through 5 above along with the bases provided in Appendix \ref{AppA} and the code in \cite{code} complete the proof. \hfill $\blacksquare$
	
	\vspace{.3in}

	\subsection*{Step 6}\label{Step 6}
	We now turn our attention to Theorem \ref{RiceqcT}.
	\\
	\\We want to find  the solutions to $ric = cT$, and we have shown that this amounts to a description of $(t_1, t_2, t_3)$ such that there exists a $c$ where $(r_1, r_2, r_3) = c(t_1, t_2, t_3)$ for some $(x,y,z)$ with $x,y > 0$ an $xy - z^2 > 0$. In the end, it is nicest to not just have a description of the $(t_1, t_2, t_3)$ without reference to the metric described by $(x,y,z)$, but also a description of the $c$ without reference to the $(x,y,z)$.
	\\
	\\As a reminder, we provide our $r_1, r_2, r_3$ without using any scale invariance below:
	\begin{align*}
		r_1 &= \frac{9 x^2 y^2-18 x y z^2+y^4+6 y^2 z^2+18 z^4}{24 \left(z^2-x y\right)^2} \\
		r_2 &= \frac{-3 x^2 \left(4 y^2-3 z^2\right)-2 x \left(y^3-12 y z^2\right)+3 y^2 z^2-6 z^4}{24 \left(z^2-x y\right)^2} \\
		r_3 &= -\frac{y z (3 x+y)^2}{24 \left(z^2-x y\right)^2}.
	\end{align*}
	Our approach to solving this is as follows. First, using the function \textcolor{red}{\textbf{\textit{FunctionRange}}} in Mathematica we determined that $r_1 > \frac{3}{8} > 0$. Using this, $(r_1, r_2, r_3) = c(t_1, t_2, t_3)$ is true if and only if $(1, \frac{r_2}{r_1}, \frac{r_3}{r_1}) = (1, \frac{t_2}{t_1}, \frac{t_3}{t_1})$. Thus, we get a description of the $(t_1, t_2, t_3)$ desired if we can describe $R = \{(1, l, m) : (1,l, m) = (1, \frac{r_2}{r_1}, \frac{r_3}{r_1})\}$. That is, we know that $(t_1, t_2, t_3)$ is a solution to $(r_1, r_2, r_3) = c(t_1, t_2, t_3)$ if and only if $(t_1, t_2, t_3)$ is such that $(1, \frac{t_2}{t_1}, \frac{t_3}{t_1}) \in R$. In addition to finding a description of $R$, we would like to find the $c$ such that $(r_1, r_2, r_3) = c(t_1, t_2, t_3)$.
	\\
	\\The subsequent material is according to the following program:
	\begin{itemize}
		\item[(a)] First, we show how one can find the $c$ value when provided a solution $(t_1, t_2, t_3)$ such that $ric = cT$. Like above, we do that in two steps: \textbf{Case 2} in which $z \not= 0$ and then \textbf{Case 1} when $z = 0$.
		\item[(b)] Our second main step is to explain how we determined the region $R$ using Mathematica and then to provide the description of $R$ with the polynomials that parameterize the regions composing $R$. This will provide us with the solution to \textbf{Case 2} of Theorem \ref{RiceqcT}.
		\item[(c)]  Lastly, we provide the solution to \textbf{Case 1} that was found using the same methodology (but greatly simplified due to fewer terms) as was used for \textbf{Case 2}.
	\end{itemize}  
	Let us begin with (a), the description of the $c$ values in our two cases for Theorem \ref{RiceqcT}.
	\\
	\\\textbf{Claim:} For $(r_1, r_2, r_3) = c(t_1, t_2, t_3)$, $c = \frac{c_0}{t_1}$ where $c_0$ is given by the implicit equation $c_0= f_1(c_0\frac{t_2}{t_1}, c_0\frac{t_3}{t_1})$ where $f_1$ is defined by Eqn.\ref{ricTsolutionso17}, unless $(t_1, t_2, t_3)$ is a multiple of $(\frac{3}{4}, r_2, r_3)$ with $r_2 > \frac{-3}{4}$ and $|r_3| = \frac{1}{2}\sqrt{\frac{3}{2}}\sqrt{4r_2 + 3}$. In this case, $c = \frac{c_1}{t_1}$ where $c_1 = \frac{3}{4}$. Moreover, our $c$ values so described are real so long as $(t_1, t_2, t_3)$ are solutions to $ric = cT$.
	\\
	\\\underline{Proof of \textbf{Claim}:}
	For a solution $(t_1, t_2, t_3)$, $(1, \frac{t_2}{t_1}, \frac{t_3}{t_1}) = (1, l, m)$ from $R$, so $\frac{1}{t_1}(t_1, t_2, t_3) = (1, l, m)$. Letting $c = \frac{1}{t_1}$ then gets us to a corresponding element of $R$, but not the $(r_1, r_2, r_3)$ in the image of $ric$. The corresponding $(r_1, r_2, r_3)$ is one that satisfies the equation $(1, \frac{r_2}{r_1}, \frac{r_3}{r_1}) = (1, l, m) = \frac{1}{t_1}(t_1, t_2, t_3)$. Observe that we then could have $(r_1, r_2, r_3) = \frac{r_1}{t_1}(t_1, t_2, t_3)$, but that would cause our $c = \frac{r_1}{t_1}$, meaning our $c$ would be dependent upon our $(x,y,z)$. A scenario which we wish to avoid. Allow us to analyze this equality more:
	\begin{align*}
		(r_1, r_2, r_3) &= \frac{r_1}{t_1}(t_1, t_2, t_3) \\
		&= (r_1, \frac{r_1t_2}{t_1}, \frac{r_1t_3}{t_1}) \\
		&\text{ which implies that } \\
		r_2 &= \frac{r_1t_2}{t_1} \\
		r_3 &= \frac{r_1t_3}{t_1}.
	\end{align*}
	Now, from Step 5, we know that if $(r_1, r_2, r_3)$ is a point on the image of $ric$, then if we have $r_2, r_3$, then we can determine $r_1$ by $r_1 = f_1(r_2, r_3) = f_1(\frac{r_1t_2}{t_1}, \frac{r_1t_3}{t_1})$ unless $r_1 = \frac{3}{4}$, $r_2 > \frac{-3}{4}$, and $|r_3| = \frac{1}{2}\sqrt{\frac{3}{2}}\sqrt{4r_2 + 3}$. We save this exception for after the rest of the solutions. 
	\\
	\\Since $(t_1, t_2, t_3)$ is provided as a solution, we can thus understand $r_1 = f_1(\frac{r_1t_2}{t_1}, \frac{r_1t_3}{t_1})$ as an implicit equation describing $r_1$ in terms of a solution $(t_1, t_2, t_3)$. Thus, we have our $c$ value in terms of $(t_1, t_2, t_3)$, namely, $\frac{c_0}{t_1}$ where $c_0$ is a solution in terms of $(t_1, t_2, t_3)$ determined by solving for $c_0$ in the implicit equation $c_0 = f_1(\frac{c_0t_2}{t_1}, \frac{c_0t_3}{t_1})$.
	\\
	\\Now for the exceptional situation. In this case, following the above procedure, we could consider $r_1 = g_1(r_2, r_3) = \frac{3}{4}$ with $c$ values given by $c = \frac{c_1}{t_1}$ where $c_1$ is a solution to $c_1 = g_1(\frac{c_1t_2}{t_1}, \frac{c_1t_3}{t_1}) = \frac{3}{4}$.
	\\
	\\To see that the $c$ values are real, consider that we were working with $(t_1, t_2, t_3)$ were solutions to $ric = cT$, implying that for a $(t_1, t_2, t_3)$ that is a solution, there is a real $c$ with the description provided. This ends the proof our \textbf{Claim}.
	\\
	\\Following the above procedure, we can see that in \textbf{Case 1} of Theorem \ref{RiceqcT}, when $z = r_3 = t_3 = 0$, we can obtain our $c$ value in a similar way. First, recall that the image of $ric$ in this case is described by $r_1 = \frac{3}{8}(5 + 16r_2 + 16r_2^2)$. Now, if $(r_1, r_2) = c(t_1, t_2)$ then, for a solution $(t_1, t_2)$ we can get $c = \frac{c_0}{t_1}$ where $c_0$ is a solution in terms of $(t_1, t_2)$ determined by solving for $c_0$ in the implicit equation $c_0 = \frac{3}{8}(5 + 16(\frac{c_0t_2}{t_1}) + 16(\frac{c_0t_2}{t_1})^2)$. This concludes the description of our $c$ values when given a solution $(t_1, t_2, t_3)$ to $ric = cT$.
	\\
	\\Now, we move to part (b) as outlined above. As a reminder, here, we discuss the region $R$, which describes the $(t_1, t_2, t_3)$ with a solution to $ric = cT$, and how we obtained said region. 
	\\
	\\
	As in Step 5, we can restrict ourselves to the setting in which $z> 0$ since we are looking for $c> 0$, so we need only consider when $r_3 > 0$. As was previously mentioned, we also provide our solution in Mathematica to the case in which $z = r_3 = 0$ (see \cite{code}).
	\\
	\\Just as before in Step 5, we us the \textcolor{red}{\textbf{\textit{Resolve}}} and \textcolor{red}{\textbf{\textit{Exists}}} functions in Mathematica to determine the image of the map described by $(1, \frac{r_2}{r_1}, \frac{r_3}{r_1}) = (1, l, m)$.
	\\
	\\The output that Mathematica provided was, as in Step 5, capable of being described more simply and is described as a collection of regions described in a piece-wise fashion. However, the description of the region $R$ is still quite messy, being described by the union of 14 different regions involving 14 different functions. These functions are, again, described implicitly as the roots of polynomials, but this time only a couple were capable of being expressed explicitly. We provide the regions below and then the polynomials needed, with the correct root being specified (note that the nth root, or Root n, is the nth root in the set of roots placed in increasing order).
	\\
	\\
	The region $R$ are those $(1, l , m)$ described by the following regions:
	\\
	$$\text{Region 1: }\begin{cases}
		m>0 & \frac{1}{2} \left(m^2-2\right)\leq l<m^2 \\
	\end{cases}
	$$
	
	$$\text{Region 2: }\begin{cases}
		0<m<\frac{1}{\sqrt{3}} & f_{275m1} \leq l<\frac{1}{3} \left(3 m^2-4\right) \\
	\end{cases}
	$$
	
	$$\text{Region 3: }\begin{cases}
		0<m\leq \sqrt{\frac{2}{3}} & f_{2213m1} \leq l\leq \frac{1}{2} \left(m^2-2\right) \\
		m>\sqrt{\frac{2}{3}} & f_{2213m1} \leq l<\frac{1}{3} \left(3 m^2-4\right) \\
	\end{cases}$$
	
	$$\text{Region 4: }\begin{cases}
		0<m\leq .625... & l = f_{4507m1} \\
		\frac{1}{\sqrt{3}}<m< .625... & l = f_{4507m2} \\
	\end{cases}$$
	
	$$\text{Region 5: }\begin{cases}
		m >\frac{1}{\sqrt{3}} & f_{275m1} \leq l<m^2 \\
	\end{cases}$$
	
	$$\text{Region 6: }\begin{cases}
		\frac{1}{\sqrt{3}}<m\leq .986... & f_{1168m1} <l<\frac{1}{3} \left(3 m^2-4\right) \\
		\frac{1}{\sqrt{3}}<m\leq  1.11... & f_{1168m2} < l < f_{25m1} \\
		1.11... < m < 1.40... & f_{2213m1} \leq l< f_{25m1} \\
		m \geq 2.22... & f_{25m1}  < l <\frac{1}{3} \left(3 m^2-4\right) \\
	\end{cases}$$
	
	$$\text{Region 7: }\begin{cases}
		\frac{1}{\sqrt{3}}<m \leq 1.11... & f_{1168m2} < l \leq f_{275m1} \\
		1.11... < m \leq 1.56... & f_{2213m1}\leq l\leq f_{275m1} \\
	\end{cases}$$
	
	$$\text{Region 8: }\begin{cases}
		\frac{1}{\sqrt{3}}<m \leq 2.22... &f_{25m1} < l \leq f_{275m1} \\
		m > 2.22...  & f_{2213m1} 	\leq l\leq f_{275m1} \\
	\end{cases}$$
	
	$$
	\text{Region 9: }\begin{cases} 
		.281... < m \leq \frac{1}{\sqrt{3}} & f_{166360m2} <l<m^2 \\
		\frac{1}{\sqrt{3}}<m\leq\sqrt{\frac{2}{3}} & f_{166360m1} <l<\frac{1}{3} \left(3 m^2-4\right) \\
		\frac{1}{\sqrt{3}}<m\leq\sqrt{\frac{2}{3}} & f_{166360m2} <l< f_{25m1} \\
		\sqrt{\frac{2}{3}}<m\leq.875... & f_{166360m1} < l < f_{25m1} \\
	\end{cases}$$

	$$\text{Region 10: }\begin{cases}
		0<m\leq .372... & f_{4507m1} <l<m^2 \\
		.372... <m\leq .556... & f_{4507m1} <l< f_{150m2} \\
		.556... <m < 1.52... & f_{2213m1} \leq l< f_{150m2} \\
		1.52... \leq m< 1.52 ... & f_{150m1} < l < f_{150m2}
	\end{cases}$$
	$$\text{Region 11: }\begin{cases}
		.372... < m\leq \frac{1}{\sqrt{3}} & f_{150m2} < l<m^2 \\
		\frac{1}{\sqrt{3}}<m\leq 1.17... & f_{150m2} < l < f_{25m1} \\
		1.17... < m < 1.40... & f_{150m1} < l < f_{25m1} \\
		1.40... \leq m \leq 1.52... & f_{150m1} < l<\frac{1}{3} \left(3 m^2-4\right) \\
	\end{cases}$$

	$$\text{Region 12: }\begin{cases}
		\frac{3 \sqrt{3}}{5}<m\leq  1.17... & f_{150m2} < l \leq f_{275m1} \\
		1.17... < m \leq 1.52... & f_{150m1} <  l  \leq f_{275m1} \\
	\end{cases}
	$$
	$$\text{Region 13: }\begin{cases}
		\frac{3 \sqrt{3}}{5}<m\leq 1.13... & f_{1168m2} < l < f_{150m2} \\
		1.09... < m \leq 1.13... & f_{2213m1} \leq l < f_{150m1} \\
	\end{cases}$$
	$$\text{Region 14: }\begin{cases}
		0< m< .423... & f_{2213m1} \leq l < f_{4507m1} \\
		.423... \leq m < .556... & f_{2213m1} \leq l \leq f_{275m1} \\
		.556... \leq m < \frac{1}{\sqrt{3}} & f_{4507m1} <l\leq f_{275m1} \\
		\frac{1}{\sqrt{3}} \leq m<.625... & f_{4507m1} < l < f_{4507m2} \\
	\end{cases}$$

	In the above regions, the numerical values bounding $m$ are approximations provided by Mathematica. For the exact expression, we reference our Mathematica code which can be found in \cite{code}. Note that the exact expressions are given as a roots of polynomials. The bounds for $l$ above can be found below and are the indicated roots to the given polynomials in $l$ where $m$ is treated as a constant.
	\\
	\\
	We provide here the polynomials in $l$ that describe the 14 regions above:
	\begin{align}
		f_{27m1} &= \text{Root 1 of the following polynomial in }l \\
		&  -2 - 75 m^2 + (6 - 180 m^2) l + (180 - 378 m^2) l^2 + (756 - 324 m^2) l^3           \notag\\
		&+(1134 - 243 m^2) l^4 + 486 l^5           \notag  \\
		\notag \\
		f_{2213m1} &= \text{Root 1 of the following polynomial in }l \notag\\
		&486 l^5+l^4 \left(1134-243 m^2\right)+l^3 \left(972 m^2+756\right)+l^2 \left(-648 m^4-1242 m^2+180\right)              \notag\\
		&+l \left(5616 m^4-3780 m^2+6\right)-2160 m^6+5112 m^4+213 m^2-2            \notag\\     
		\notag \\
		f_{4507m1} &= \text{Root 1 of the following polynomial in }l \notag\\
		&4 - 507 m^2 + (51 - 1404 m^2) l + (252 - 1674 m^2) l^2 + (594 - 972 m^2) l^3          \notag\\
		&+(648 - 243 m^2) l^4 + 243 l^5            \notag\\
		\notag \\
		f_{4507m2} &= \text{Root 2 of the following polynomial in }l \notag\\
		&4 - 507 m^2 + (51 - 1404 m^2) l + (252 - 1674 m^2) l^2 + (594 - 972 m^2) l^3           \notag\\
		&+(648 - 243 m^2) l^4 + 243 l^5           \notag\\
		\notag \\
		f_{25m1} &= \text{Root 1 of the following polynomial in }l \notag\\
		&-25 m^2 + (1 + 60 m^2) l + (12 - 126 m^2) l^2 + (54 + 108 m^2) l^3 + (108 - 81 m^2) l^4 + 81 l^5           \notag\\
		\notag \\
		f_{1168m1} &= \text{Root 1 of the following polynomial in }l \notag\\
		&1 - 168 m^2 + 144 m^4 + (12 + 144 m^2) l + (54 + 216 m^2) l^2
		+108 l^3 + 81 l^4         \notag\\
		\notag \\
		f_{1168m2} &=  \text{Root 2 of the following polynomial in }l \notag\\
		&1 - 168 m^2 + 144 m^4 + (12 + 144 m^2) l + (54 + 216 m^2) l^2 + 
		108 l^3 + 81 l^4          \notag\\
		\notag
	\end{align}

	\begin{align*}
		f_{1166360m1} &= \text{Root 1 of the following polynomial in }l \notag\\
		&16 - 6360 m^2 + 47961 m^4 + (216 + 8478 m^2 + 149796 m^4) l          \notag\\
		&+(1161 + 54432 m^2 + 176094 m^4) l^2 + (3132 + 64476 m^2 + 92340 m^4) l^3            \notag\\
		&+(4374 + 29160 m^2 + 18225 m^4) l^4 +(2916 + 4374 m^2) l^5 + 729 l^6           \notag\\
		\notag \\
		f_{166360m2} &= \text{Root 2 of the following polynomial in }l \notag\\
		&16 - 6360 m^2 + 47961 m^4 + (216 + 8478 m^2 + 149796 m^4) l          \notag\\
		&+(1161 + 54432 m^2 + 176094 m^4) l^2 + (3132 + 64476 m^2 + 92340 m^4) l^3            \notag\\
		&+(4374 + 29160 m^2 + 18225 m^4) l^4 +(2916 + 4374 m^2) l^5 + 729 l^6           \notag\\
		\notag \\
		f_{150m1} &= \text{Root 1 of the following polynomial in }l \notag\\
		&-150 m^2 + 75 m^4 + (2 - 231 m^2 + 540 m^4) l + (42 + 4806 m^2 + 3402 m^4) l^2          \notag\\
		&+(378 + 14742 m^2 + 8748 m^4) l^3 +(1890 - 6966 m^2 + 19683 m^4) l^4          \notag\\
		&+(5670 - 45927 m^2) l^5 + (10206 - 13122 m^2) l^6 + 10206 l^7 + 4374 l^8 &           \notag\\
		\notag \\
		f_{150m2} &= \text{Root 2 of the following polynomial in }l \notag\\
		&-150 m^2 + 75 m^4 + (2 - 231 m^2 + 540 m^4) l + (42 + 4806 m^2 + 3402 m^4) l^2          \notag\\
		&+(378 + 14742 m^2 + 8748 m^4) l^3 +(1890 - 6966 m^2 + 19683 m^4) l^4          \notag\\
		&+(5670 - 45927 m^2) l^5 + (10206 - 13122 m^2) l^6 + 10206 l^7 + 4374 l^8 &           \notag
		\\\notag
	\end{align*}
	Of the above roots, there are only two for which Mathematica could produce a radical description:
	$$f_{1168m1} = \frac{1}{3} \left(-2 \sqrt{-3 m^2-2 \sqrt{3} m}-1\right)$$
	$$f_{1168m2} = \frac{1}{3} \left(2 \sqrt{-3 m^2-2 \sqrt{3} m}-1\right).$$
	Having the above regions and functions, we have the solution to $ric = cT$ for $z, t_3 > 0$. To find solutions for $t_3\not= 0$, as before, one must exchange conditions on $t_3$ to be conditions on $|t_3|$, meaning that for $t_3 \not= 0$, we want $(t_1, t_2, t_3)$ to be such that $\frac{t_2}{t_1} = l $ and $\frac{|t_3|}{t_1} = m$ in the 14 regions above. Thus, we have solution to \textbf{Case 2} of Theorem \ref{RiceqcT}.
	\\
	\\Lastly, we discuss part (c) as was outlined above. Here, we turn our attention to the description of $(t_1, t_2, t_3)$ in \textbf{Case 1} of Theorem \ref{RiceqcT}.
	\\
	\\
	As mentioned before, we provide the Mathematica code for \textbf{Case 1} of Theorem \ref{RiceqcT} in \cite{code}. The proof for this case follows the same program as \textbf{Case 2}, but with $z = r_3 = 0$ which greatly simplifies things. Indeed, if one looks to our Mathematica work in the $z = 0$ setting, one can notice that finding the description of $R = \{(1, l): (1, l) = (1, \frac{r_2}{r_1})\}$ using \textcolor{red}{\textbf{\textit{Exists}}} and \textcolor{red}{\textbf{\textit{Resolve}}} amounts to one simple region described by:
	$$\frac{1}{3} \left(-\sqrt{5}-2\right)\leq l<0.$$
	Therefore, $(t_1, t_2)$ providing solutions in the case of $z = 0$ are those which satisfy $$\frac{1}{3} \left(-\sqrt{5}-2\right)\leq \frac{t_2}{t_1} <0.$$
	This provides us with the desired solution of \textbf{Case 1} in Theorem \ref{RiceqcT}.
	\\
	\\Having solutions provided above and the proof of our \textbf{Claim} above, we have a description of our $(t_1, t_2, t_3)$ and $c>0$ such that $(r_1, r_2, r_3) = c(t_1, t_2, t_3)$ for some $(x,y,z)$ with $xy - z^2 > 0$ and $x, y > 0$.
	\\
	\\\underline{Proof of Theorem \ref{RiceqcT}}: Steps 1 through 6 above along with the work provided in Mathematica in Appendix \ref{AppA} and the code in \cite{code} complete the proof. \hfill $\blacksquare$
	\begin{remark}
		In the following pages, we provide images depicting $R$. The images indicate that there is still some simplification that could happen. We remind the reader that $R$ describes the solutions for $z > 0$ which in turn gives us solutions for $z \not= 0$. The $z < 0$ solutions would just be reflections of the graphs seen about the horizontal axes.
		\\
		\\
	\end{remark}
	
	\noindent In the following images, we depict the region $R$ with the vertical axis being the $m$ axis and the horizontal axis the $l$ axis. The first image is all of the regions shown at the same time. We remove various regions to show how there appears to be great overlap between regions, indicating that further simplification of the description of $R$ could be made.
	\\
	
	\includegraphics[width=\textwidth]{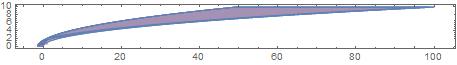}
	\vfill
	
	\noindent Below we have $R$ with all of the regions except for the second part of Region 3 and Region 5. The larger portion that is orange is Region 1.
	\\

	\includegraphics[width=\textwidth]{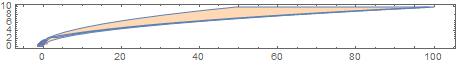}
	\vfill

	\noindent Below we have $R$ with all of the regions except for Region 1 and Region 5. The larger portion in blue is the second part of Region 3.
	\\

	\includegraphics[width=\textwidth]{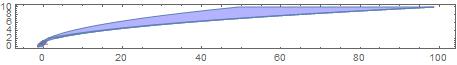}
	\vfill

	\noindent Below we have $R$ with all of the regions except for Region 1 and the second part of Region 3. The larger grey part is Region 5.
	\\

	\includegraphics[width=\textwidth]{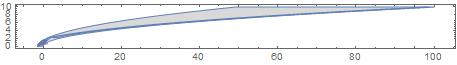}
	\vfill
	
	\noindent Below we have $R$ with all of the regions except for Region 1, the second part of Region 3, and Region 5. The thin part going up with greater $m$ values is the second part of Region 8, and the thin part going up with lesser $m$ values is the last region in Region 6.
	\\

	\includegraphics[width=\textwidth]{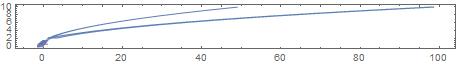}
	\vfill
	
	\noindent Below we have $R$ with all of the regions except for Region 1, the second part of Region 3, Region 5, and the second part of Region 8. The thin part going up is the last region in Region 6.
	\\

	\includegraphics[width=\textwidth]{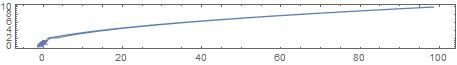}
	\vfill

	\pagebreak
	
	\noindent Below we have $R$ with all of the regions except for Region 1, the second part of Region 3, Region 5, the last part of Region 6, and the second part of Region 8. Notice by observing the bounds that the amount of $R$ seen has decreased significantly, but there is still much overlap between regions.
	\\
	
	\includegraphics[width=\textwidth]{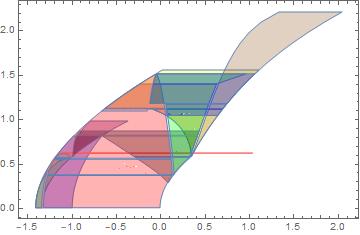}
	
	\newpage

	\appendix

	\section{Appendix A}\label{AppA}
	\vspace{.2in}
	
	\noindent Here we provide the basis of $\p_1$ and $\p_2$ in $\so(1,7) = \g_2 \oplus \p_1 \oplus \p_2$ that we ultimately use to compute $ric$. This basis is orthonormal with respect to our fixed inner product $\langle . , . \rangle$.
	\\
	\\
	\\
	Our basis for $\p_1$:
	\\
	\\1
	$\left[\begin{matrix}0 & 0 & 0 & 0 & 0 & 0 & 0 & 0\\0 & 0 & 0 & 0 & 0 & 0 & 0 & 0\\0 & 0 & 0 & \frac{1}{6} & 0 & 0 & 0 & 0\\0 & 0 & - \frac{1}{6} & 0 & 0 & 0 & 0 & 0\\0 & 0 & 0 & 0 & 0 & \frac{1}{6} & 0 & 0\\0 & 0 & 0 & 0 & - \frac{1}{6} & 0 & 0 & 0\\0 & 0 & 0 & 0 & 0 & 0 & 0 & - \frac{1}{6}\\0 & 0 & 0 & 0 & 0 & 0 & \frac{1}{6} & 0\end{matrix}\right]$
	\hfill
	2
	$\left[\begin{matrix}0 & 0 & 0 & 0 & 0 & 0 & 0 & 0\\0 & 0 & 0 & - \frac{1}{6} & 0 & 0 & 0 & 0\\0 & 0 & 0 & 0 & 0 & 0 & 0 & 0\\0 & \frac{1}{6} & 0 & 0 & 0 & 0 & 0 & 0\\0 & 0 & 0 & 0 & 0 & 0 & \frac{1}{6} & 0\\0 & 0 & 0 & 0 & 0 & 0 & 0 & \frac{1}{6}\\0 & 0 & 0 & 0 & - \frac{1}{6} & 0 & 0 & 0\\0 & 0 & 0 & 0 & 0 & - \frac{1}{6} & 0 & 0\end{matrix}\right]
	$\hfill
	\\
	3
	$\left[\begin{matrix}0 & 0 & 0 & 0 & 0 & 0 & 0 & 0\\0 & 0 & \frac{1}{6} & 0 & 0 & 0 & 0 & 0\\0 & - \frac{1}{6} & 0 & 0 & 0 & 0 & 0 & 0\\0 & 0 & 0 & 0 & 0 & 0 & 0 & 0\\0 & 0 & 0 & 0 & 0 & 0 & 0 & \frac{1}{6}\\0 & 0 & 0 & 0 & 0 & 0 & - \frac{1}{6} & 0\\0 & 0 & 0 & 0 & 0 & \frac{1}{6} & 0 & 0\\0 & 0 & 0 & 0 & - \frac{1}{6} & 0 & 0 & 0\end{matrix}\right]$ \hfill
	4
	$\left[\begin{matrix}0 & 0 & 0 & 0 & 0 & 0 & 0 & 0\\0 & 0 & 0 & 0 & 0 & - \frac{1}{6} & 0 & 0\\0 & 0 & 0 & 0 & 0 & 0 & - \frac{1}{6} & 0\\0 & 0 & 0 & 0 & 0 & 0 & 0 & - \frac{1}{6}\\0 & 0 & 0 & 0 & 0 & 0 & 0 & 0\\0 & \frac{1}{6} & 0 & 0 & 0 & 0 & 0 & 0\\0 & 0 & \frac{1}{6} & 0 & 0 & 0 & 0 & 0\\0 & 0 & 0 & \frac{1}{6} & 0 & 0 & 0 & 0\end{matrix}\right]$\hfill
	\\
	5
	$\left[\begin{matrix}0 & 0 & 0 & 0 & 0 & 0 & 0 & 0\\0 & 0 & 0 & 0 & \frac{1}{6} & 0 & 0 & 0\\0 & 0 & 0 & 0 & 0 & 0 & 0 & - \frac{1}{6}\\0 & 0 & 0 & 0 & 0 & 0 & \frac{1}{6} & 0\\0 & - \frac{1}{6} & 0 & 0 & 0 & 0 & 0 & 0\\0 & 0 & 0 & 0 & 0 & 0 & 0 & 0\\0 & 0 & 0 & - \frac{1}{6} & 0 & 0 & 0 & 0\\0 & 0 & \frac{1}{6} & 0 & 0 & 0 & 0 & 0\end{matrix}\right]$ \hfill
	6
	$\left[\begin{matrix}0 & 0 & 0 & 0 & 0 & 0 & 0 & 0\\0 & 0 & 0 & 0 & 0 & 0 & 0 & \frac{1}{6}\\0 & 0 & 0 & 0 & \frac{1}{6} & 0 & 0 & 0\\0 & 0 & 0 & 0 & 0 & - \frac{1}{6} & 0 & 0\\0 & 0 & - \frac{1}{6} & 0 & 0 & 0 & 0 & 0\\0 & 0 & 0 & \frac{1}{6} & 0 & 0 & 0 & 0\\0 & 0 & 0 & 0 & 0 & 0 & 0 & 0\\0 & - \frac{1}{6} & 0 & 0 & 0 & 0 & 0 & 0\end{matrix}\right]$ \hfill
	\\
	7
	$\left[\begin{matrix}0 & 0 & 0 & 0 & 0 & 0 & 0 & 0\\0 & 0 & 0 & 0 & 0 & 0 & - \frac{1}{6} & 0\\0 & 0 & 0 & 0 & 0 & \frac{1}{6} & 0 & 0\\0 & 0 & 0 & 0 & \frac{1}{6} & 0 & 0 & 0\\0 & 0 & 0 & - \frac{1}{6} & 0 & 0 & 0 & 0\\0 & 0 & - \frac{1}{6} & 0 & 0 & 0 & 0 & 0\\0 & \frac{1}{6} & 0 & 0 & 0 & 0 & 0 & 0\\0 & 0 & 0 & 0 & 0 & 0 & 0 & 0\end{matrix}\right]$
	\\
	\\
	\\
	\\
	\\
	\\
	\\
	\\Our basis for $\p_2$
	\\
	\\8
	$\left[\begin{matrix}0 & \frac{\sqrt{3}}{6} & 0 & 0 & 0 & 0 & 0 & 0\\\frac{\sqrt{3}}{6} & 0 & 0 & 0 & 0 & 0 & 0 & 0\\0 & 0 & 0 & 0 & 0 & 0 & 0 & 0\\0 & 0 & 0 & 0 & 0 & 0 & 0 & 0\\0 & 0 & 0 & 0 & 0 & 0 & 0 & 0\\0 & 0 & 0 & 0 & 0 & 0 & 0 & 0\\0 & 0 & 0 & 0 & 0 & 0 & 0 & 0\\0 & 0 & 0 & 0 & 0 & 0 & 0 & 0\end{matrix}\right]$ \hfill
	9
	$\left[\begin{matrix}0 & 0 & \frac{\sqrt{3}}{6} & 0 & 0 & 0 & 0 & 0\\0 & 0 & 0 & 0 & 0 & 0 & 0 & 0\\\frac{\sqrt{3}}{6} & 0 & 0 & 0 & 0 & 0 & 0 & 0\\0 & 0 & 0 & 0 & 0 & 0 & 0 & 0\\0 & 0 & 0 & 0 & 0 & 0 & 0 & 0\\0 & 0 & 0 & 0 & 0 & 0 & 0 & 0\\0 & 0 & 0 & 0 & 0 & 0 & 0 & 0\\0 & 0 & 0 & 0 & 0 & 0 & 0 & 0\end{matrix}\right]$ \hfill
	\\
	10
	$\left[\begin{matrix}0 & 0 & 0 & \frac{\sqrt{3}}{6} & 0 & 0 & 0 & 0\\0 & 0 & 0 & 0 & 0 & 0 & 0 & 0\\0 & 0 & 0 & 0 & 0 & 0 & 0 & 0\\\frac{\sqrt{3}}{6} & 0 & 0 & 0 & 0 & 0 & 0 & 0\\0 & 0 & 0 & 0 & 0 & 0 & 0 & 0\\0 & 0 & 0 & 0 & 0 & 0 & 0 & 0\\0 & 0 & 0 & 0 & 0 & 0 & 0 & 0\\0 & 0 & 0 & 0 & 0 & 0 & 0 & 0\end{matrix}\right]$ \hfill
	11
	$\left[\begin{matrix}0 & 0 & 0 & 0 & \frac{\sqrt{3}}{6} & 0 & 0 & 0\\0 & 0 & 0 & 0 & 0 & 0 & 0 & 0\\0 & 0 & 0 & 0 & 0 & 0 & 0 & 0\\0 & 0 & 0 & 0 & 0 & 0 & 0 & 0\\\frac{\sqrt{3}}{6} & 0 & 0 & 0 & 0 & 0 & 0 & 0\\0 & 0 & 0 & 0 & 0 & 0 & 0 & 0\\0 & 0 & 0 & 0 & 0 & 0 & 0 & 0\\0 & 0 & 0 & 0 & 0 & 0 & 0 & 0\end{matrix}\right]$ \hfill
	\\
	12
	$\left[\begin{matrix}0 & 0 & 0 & 0 & 0 & \frac{\sqrt{3}}{6} & 0 & 0\\0 & 0 & 0 & 0 & 0 & 0 & 0 & 0\\0 & 0 & 0 & 0 & 0 & 0 & 0 & 0\\0 & 0 & 0 & 0 & 0 & 0 & 0 & 0\\0 & 0 & 0 & 0 & 0 & 0 & 0 & 0\\\frac{\sqrt{3}}{6} & 0 & 0 & 0 & 0 & 0 & 0 & 0\\0 & 0 & 0 & 0 & 0 & 0 & 0 & 0\\0 & 0 & 0 & 0 & 0 & 0 & 0 & 0\end{matrix}\right]$ \hfill
	13
	$\left[\begin{matrix}0 & 0 & 0 & 0 & 0 & 0 & \frac{\sqrt{3}}{6} & 0\\0 & 0 & 0 & 0 & 0 & 0 & 0 & 0\\0 & 0 & 0 & 0 & 0 & 0 & 0 & 0\\0 & 0 & 0 & 0 & 0 & 0 & 0 & 0\\0 & 0 & 0 & 0 & 0 & 0 & 0 & 0\\0 & 0 & 0 & 0 & 0 & 0 & 0 & 0\\\frac{\sqrt{3}}{6} & 0 & 0 & 0 & 0 & 0 & 0 & 0\\0 & 0 & 0 & 0 & 0 & 0 & 0 & 0\end{matrix}\right]$ \hfill
	\\
	14
	$\left[\begin{matrix}0 & 0 & 0 & 0 & 0 & 0 & 0 & \frac{\sqrt{3}}{6}\\0 & 0 & 0 & 0 & 0 & 0 & 0 & 0\\0 & 0 & 0 & 0 & 0 & 0 & 0 & 0\\0 & 0 & 0 & 0 & 0 & 0 & 0 & 0\\0 & 0 & 0 & 0 & 0 & 0 & 0 & 0\\0 & 0 & 0 & 0 & 0 & 0 & 0 & 0\\0 & 0 & 0 & 0 & 0 & 0 & 0 & 0\\\frac{\sqrt{3}}{6} & 0 & 0 & 0 & 0 & 0 & 0 & 0\end{matrix}\right]$
	\\

	\newpage

	\bibliographystyle{alpha}
	
	\bibliography{ONTHEPRESCRIBEDRICCICURVATUREOFNONCOMPACT_ArXiv1}

\end{document}